\documentclass[a4paper]{article}
\usepackage[margin=25mm]{geometry}
\usepackage[version=4]{mhchem}
\usepackage{bm}
\usepackage{xifthen}
\usepackage{graphicx}
\usepackage{caption}
\usepackage{subcaption}
\usepackage{float}
\usepackage[hidelinks]{hyperref}
\usepackage{amsfonts}
\usepackage{verbatim,xcolor}

\hypersetup{
    breaklinks=true,
    colorlinks=true, 
    linkcolor=blue,  
    citecolor=blue,    
    filecolor=blue,    
    urlcolor=blue          
}

\makeatletter
\def\@seccntformat#1{\@ifundefined{#1@cntformat}%
   {\csname the#1\endcsname\quad}  
   {\csname #1@cntformat\endcsname}
}
\let\oldappendix\appendix 
\renewcommand\appendix{%
    \oldappendix
    \newcommand{\section@cntformat}{\appendixname~\thesection:~}
}
\makeatother

\usepackage{fancyhdr}
\usepackage[capitalize,nameinlink]{cleveref}[0.19]
\crefname{section}{Section}{sections}
\crefname{subsection}{Subsection}{subsections}
\Crefname{section}{Section}{Sections}
\Crefname{subsection}{Subsection}{Subsections}

\Crefname{figure}{Figure}{Figures}

\crefformat{equation}{\textup{#2(#1)#3}}
\crefrangeformat{equation}{\textup{#3(#1)#4--#5(#2)#6}}
\crefmultiformat{equation}{\textup{#2(#1)#3}}{ and \textup{#2(#1)#3}}
{, \textup{#2(#1)#3}}{, and \textup{#2(#1)#3}}
\crefrangemultiformat{equation}{\textup{#3(#1)#4--#5(#2)#6}}%
{ and \textup{#3(#1)#4--#5(#2)#6}}{, \textup{#3(#1)#4--#5(#2)#6}}{, and \textup{#3(#1)#4--#5(#2)#6}}

\Crefformat{equation}{#2Equation~\textup{(#1)}#3}
\Crefrangeformat{equation}{Equations~\textup{#3(#1)#4--#5(#2)#6}}
\Crefmultiformat{equation}{Equations~\textup{#2(#1)#3}}{ and \textup{#2(#1)#3}}
{, \textup{#2(#1)#3}}{, and \textup{#2(#1)#3}}
\Crefrangemultiformat{equation}{Equations~\textup{#3(#1)#4--#5(#2)#6}}%
{ and \textup{#3(#1)#4--#5(#2)#6}}{, \textup{#3(#1)#4--#5(#2)#6}}{, and \textup{#3(#1)#4--#5(#2)#6}}

\crefdefaultlabelformat{#2\textup{#1}#3}

\usepackage[backend=biber,style=numeric-comp,giveninits=true,doi=false,url=false,isbn=false,maxbibnames=5,eprint=false,sorting=none]{biblatex} 
\newcommand{\footremember}[2]{%
    \footnote{#2}
    \newcounter{#1}
    \setcounter{#1}{\value{footnote}}%
}

\numberwithin{equation}{section}
\numberwithin{table}{section}
\numberwithin{figure}{section}

\title{Accurate stochastic simulation of nonlinear reactions between closest particles\footnote{Funding: This work has been partially supported by the Australian Research Council through the Future Fellowship grant FT220100496 and Discovery Project grant DP22010316.}}

\author{%
  Taylor Kearney\footremember{alley}{School of Mathematics, Monash University, 9 Rainforest walk, 3800 Clayton, Victoria, Australia, Taylor.Kearney1@monash.edu}%
   \and Ricardo Ruiz-Baier\footremember{trailer}{School of Mathematics, Monash University, 9 Rainforest walk, 3800 Clayton, Victoria, Australia, Ricardo.RuizBaier@monash.edu}
  \and Mark B. Flegg\footremember{underbridge}{School of Mathematics, Monash University, 9 Rainforest walk, 3800 Clayton, Victoria, Australia,  Mark.Flegg@monash.edu}
 }
\date{}

\providecommand{\keywords}[1]
{
  \small	
  \textbf{{Keywords---}} #1
}
\newcommand{\pos}[2]{\boldsymbol{#1}_{#2}}
\newcommand{\posbar}[2]{\bar{\boldsymbol{#1}}_{#2}}
\newcommand{\parentheses}[3]{\left#1#2\right#3}

\newcommand{\prob}[3][]{\ifthenelse{\isempty{#1}}{P\left(#2,#3\right)}{P\left(#1,#2,#3\right)}}

\newcommand{\probtwo}[3][]{\ifthenelse{\isempty{#1}}{P_2\left(#2,#3\right)}{P_2\left(#1,#2,#3\right)}}
\newcommand{\probGiven}[4][]{\ifthenelse{\isempty{#1}}{P\parentheses{(}{#2,#3|#4}{)}}{P\parentheses{(}{#1,#2,#3|#4}{)}}}
\newcommand{\phiFunc}{\phi\left(\pos{\eta}{2},t\right)}
\newcommand{\phiFuncPrime}{\phi\left(\pos{\eta'}{2},t\right)}
\newcommand{\PhiFunc}{\Phi\left(\pos{\eta}{2},t\right)}

\newcommand{\gFunc}{g\left(\pos{\eta}{2},t\right)}
\newcommand{\probbar}[3][]{\ifthenelse{\isempty{#1}}{\bar{P}\left(#2,#3\right)}{\bar{P}\left(#1,#2,#3\right)}}
\newcommand{\steadyProb}{P\left(r_1,r_2\right)}
\newcommand{\Q}{Q\left(r_1,r_2\right)}
\newcommand{\s}{\boldsymbol{s}\left(r_1,r_2\right)}
\newcommand{\steadyProbfe}{P_h\left(r_1,r_2\right)}

\newcommand{\sfe}{\boldsymbol{s}_h\left(r_1,r_2\right)}
\newcommand{\sigmamax}{\sigma_{\text{max}}}
\newcommand{\sigmaeff}{\sigma_{\text{E}}}
\newcommand{\gammaeff}{\Gamma_{\text{E}}}
\newcommand{\sigmainit}{\sigma_{\text{D}}}
\newcommand{\gammainit}{\Gamma_{\text{D}}}

\newcommand{\bs}{\boldsymbol{s}}
\newcommand{\btau}{\boldsymbol{\tau}}

\begin{document}

\maketitle

\begin{abstract}
We study a system of diffusing point particles in which any triplet of particles reacts and is removed from the system when the relative proximity of the constituent particles satisfies a predefined condition. Proximity-based reaction conditions of this kind are commonly used in particle-based simulations of chemical kinetics to mimic bimolecular reactions, those involving just two reactants, and have been extensively studied. The rate at which particles react within the system is determined by the reaction condition and particulate diffusion. In the bimolecular case, analytic relations exist between the reaction rate and the distance at which particles react allowing modellers to tune the rate of the reaction within their simulations by simply altering the reaction condition. However, generalising proximity-based reaction conditions to trimolecular reactions, those involving three particles, is more complicated because it requires understanding the distribution of the closest diffusing particle to a point in the vicinity of a spatially dependent absorbing boundary condition. We find that in this case the evolution of the system is described by a nonlinear partial integro-differential equation with no known analytic solution, which makes it difficult to relate the reaction rate to the reaction condition. To resolve this, we use singular perturbation theory to obtain a leading-order solution and show how to derive an approximate expression for the reaction rate. We then use finite element methods to quantify the higher-order corrections to this solution and the reaction rate, which are difficult to obtain analytically. Leveraging the insights gathered from this analysis, we demonstrate how to correct for the errors that arise from adopting the approximate expression for the reaction rate, enabling for the construction of more accurate particle-based simulations than previously possible.
\end{abstract}


\section{Introduction}
Transport-dependent phenomena are ubiquitous in many chemical systems, particularly in intracellular environments. The most fundamental example of molecular transport is the thermally driven random motion that results in a net movement from regions of high concentration to those with lower concentrations, known as diffusion \cite{Jacobs_diffusion_book}. Due to its fundamental nature, diffusion is both ubiquitous and essential in biology \cite{mogre2020getting}. It plays a central role in the operation of microorganisms \cite{Koch1990} and is vital in numerous cellular processes, including metabolism \cite{jones1986intracellular,kinsey2011molecules} and signalling \cite{kholodenko2008giving,ye2013single,stemCellSignaling}. More broadly, diffusion facilitates the mixing of reactants on small spatial scales, promoting the molecular collisions necessary for biochemical reactions to occur \cite{north1966diffusion, molecularBioCell}. If reactant diffusion is sufficiently slow, it largely determines the reaction rate, leading to a \textit{diffusion-limited} reaction. This has driven significant interest in diffusion-limited reactions and their kinetics have been extensively studied at various levels of fidelity.
\par 

Many models consider the kinetics of diffusion-limited reactions to be deterministic, assuming that the reactants can be well approximated using continuous concentrations. These deterministic models either use ordinary differential equations (ODEs) to investigate only the temporal evolution of reactant concentrations \cite{murray2007mathematical, chen2010classic} or employ partial differential equations (PDEs), known as reaction-diffusion equations, to explicitly incorporate reactant diffusion \cite{soh2010reaction, doi:10.1126/science.1179047}. Although such models have been widely studied, they provide an inadequate description of many biological systems whose kinetics are inherently stochastic \cite{stochastic_processes_bio_book, Noisy_cell_business, doi:10.1126/science.1147888,10.1371/journal.pcbi.1002010} due to their molecular origins. Moreover, molecular populations are often so low or so spatially localised that meaningful concentrations cannot be defined \cite{laws_chem_living_cells, erban2009stochastic}.
\par

These shortcomings have motivated the development of a range of models that incorporate stochastic effects as a core component. Gillespie's stochastic simulation method \cite{GILLESPIE1976403, Gillespie1977} is one of the earliest and most influential models in this class, enabling the computation of the stochastic time evolution of molecular populations in a manner consistent with the Chemical Master Equation. Gillespie's original algorithm applies only to spatially homogeneous, or well-mixed chemical systems, but it can be extended to spatially inhomogeneous systems by dividing the domain into spatially homogeneous voxels or compartments that are coupled through diffusive transfers \cite{PhysRevE.71.041103, Gillespie_spatial_inhomo}.

This approach has become commonplace in a broader class of models often referred to as reaction-diffusion master equation (RDME) models \cite{Isaacson_RDME}, which discretise space in the same manner but do not necessarily rely on Gillespie's algorithm to evolve voxel populations \cite{Early_compartment, mesoRD_og, mesoRD,doi:10.1137/080721388, hepburn2012steps}. RDME models provide a computationally efficient way to incorporate spatial dependence in stochastic models. However, the size of the voxels is inherently limited by the assumption that they are well mixed. In many cases, finer spatial resolutions are desirable to accurately capture system behaviour, leading to more computationally intensive RDME approaches \cite{CRDME_Isaacson}.

Perhaps the most natural choice for models of reaction-diffusion systems are those that aim to faithfully replicate the actual dynamics of the system by explicitly treating molecules as individuals undergoing a reaction-diffusion process. Molecular dynamics simulations \cite{molecular_dynamics_simulations} take this notion to the extreme, attempting to accurately reproduce the physical system with atomic detail. Although they are highly accurate, these simulations can be prohibitively expensive due to the sheer number of degrees of freedom that must be modelled; many of which arise from the need to explicitly track solvent molecules, account for numerous possible interactions, and capture intramolecular timescales \cite{coarse_grained_MD}.

To partially mitigate these computational challenges, particle-based models leverage simplifying assumptions to improve efficiency \cite{feig2019whole} and are widely used in the literature \cite{schoneberg2014simulation}. These models avoid the explicit treatment of solvent molecules and describe the relevant reactants as points undergoing isotropic diffusion on a continuous domain. In the absence of spatial morphology for individual molecules, these approaches typically favour idealised proximity-based (and sometimes energy-potential-based) reaction conditions that are simple to implement, making such simulations feasible to run on personal computing hardware.

This approach was largely inspired by Smoluchowski \cite{smoluchowski1917versuch}, who postulated that bimolecular reactions could be modelled by assuming that two molecules undergo a reaction if they come within a predefined distance $\sigma$, called the \textit{reaction radius}. Although simple, Smoluchowski's reaction condition has proven remarkably robust and forms the basis of many prominent software packages for particle-based simulation of reaction-diffusion systems, including Green’s function reaction dynamics (GFRD) \cite{gfrd1,gfrd2}, enhanced Green's function reaction dynamics (eGFRD) \cite{egfrd1,egfrd2}, Smoldyn \cite{Andrews_2004}, and ReaDDy \cite{schoneberg2013readdy}.

For a bimolecular reaction between two chemical species $A$ and $B$, Smoluchowski's model imposes that the rate of reaction is simply the number of successful `collisions' between molecules of $A$ and $B$ at a distance $\sigma$ per unit time. This reaction rate is inherently stochastic because it depends on the motion of individual diffusing molecules. However, if we were to average the reaction rates obtained from a large ensemble of independent molecules of $A$ surrounded by molecules of $B$, then the expected reaction rate should converge to the corresponding macroscopic rate. This provides a matching condition that can be used to determine the currently unknown reaction radius $\sigma$.\par

To derive this condition, we consider a system of extremely large (relative to molecular sizes) volume $V$ that contains a single isolated molecule of $A$ surrounded by molecules of $B$. For this system, the joint probability density $P_S\left(r_1,t\right)$ for finding a molecule of $B$, that has never collided with $A$, at a given position at time $t$ only depends on the radial distance $r_1$ between this molecule and $A$ and evolves according to the diffusion equation,
\begin{equation}\label{eq:smol_diffusion}
\frac{dP_S\left(r_1,t\right)}{dt} = \frac{\hat{D}_1}{r1^2}\frac{\partial}{\partial r_1}\left(r1^2 \frac{\partial P_S\left(r_1,t\right)}{\partial r_1}\right) ,
\end{equation}
where $\hat{D}_1 = D_0 + D_1$ is the relative diffusion coefficient between molecules $A$ and $B$ (the sum of their diffusion coefficients). Smoluchowski's reaction condition specifies that molecules of $B$ must be removed from the system at a distance $r_1 = \sigma$ from $A$ which manifests as an absorbing boundary condition
\begin{equation}\label{eq:smol_boundary_condition}
    P_S\left(r_1,t\right) = 0 \text{ on } r_1=\sigma.
\end{equation}
The expected reaction rate is given by the diffusive flux over this boundary which rapidly approaches the pseudo-equilibrium value \cite{smoluchowski1917versuch, comprehensiveChemicalKinetics}
\begin{equation}\label{eq:smol_steady_rate_constant}
    K_S\left(\sigma\right) = \frac{k_2}{V} = \frac{4\pi \sigma \hat{D}_1}{V},
\end{equation}
where $k_2$ is the bimolecular rate of reaction (first passage rate for a collision). To obtain this pseudo-equilibrium, it is enough to assume that as $r_1 \rightarrow \infty$ the probability density is constant at $V^{-1}$ since $\sigma$ is small.\par 

Equation \cref{eq:smol_steady_rate_constant} gives the reaction rate for a single molecule of $A$ per molecule of $B$ and is used to determine $\sigma$ based on a known second-order rate constant $k_2$. To recover the macroscopic reaction rate we assume that molecules of $A$ can be considered independently, so that the expected reaction rate for a system containing multiple molecules is given by Equation \cref{eq:smol_steady_rate_constant} multiplied the number of molecules of $A$ and the number of molecules of $B$. Note that often Equation \cref{eq:smol_steady_rate_constant} is instead multiplied by the concentration of $A$ (rather than the number of $A$ molecules), since for an elementary bimolecular reaction the resulting rate of change for the concentration can be directly compared to the macroscopic rate given by the law of mass action, which states the rate is given by the product of the reactant concentrations multiplied by the second-order rate constant $k_2$ associated with the reaction.\par

The validity of Smoluchowski's reaction condition, and the many derivatives it has inspired \cite{collins1949diffusion, Doi_1976, erban2009stochastic}, is contingent upon the distance at which molecules interact being sufficiently small. Specifically, for Smoluchowski's model, the reaction radius $\sigma$ must be sufficiently small when compared to the expected separation between the molecules \cite{our_first_paper}. In other words, we require the system to be sufficiently sparse or the affinity of $A$ for $B$ to be sufficiently low so that the associated reaction radius is small. \textcite{collins1949diffusion} alleviate this requirement by supposing that not every collision results in a reaction which amounts to replacing Smoluchowski's absorbing boundary in Equation \cref{eq:smol_boundary_condition} with a partially reflective condition. Similarly, \textcite{Doi_1976} allows reactants to approach arbitrarily close to each other and only requires that reactions occur at a fixed rate per unit time once reactants are separated by a distance less than the analogous reaction radius parameter. These alternative reaction conditions enable the distance at which molecules react to be increased for a fixed reaction rate when compared to Smoluchowski's reaction condition, but in each case the reaction condition is still fundamentally a proximity-dependent interaction between molecules that can be summarised by a single parameter, analogous to $\sigma$, which determines the spatial scale at which interactions occur \cite{flegg2016smoluchowski}. Moreover, for our assumptions regarding the independence of reactants of the same species to be accurate, the interaction distance must still be small regardless of the reaction condition adopted.\par

\textcite{flegg2016smoluchowski} demonstrated that proximity-based reaction conditions equivalent to Smoluchowski's can also be applied to higher-order reactions (those involving more than two reactants), while \textcite{our_first_paper} demonstrated proximity-based trimolecular reaction conditions can produce kinetics that cannot be described by the law of mass action. In this article, we consider the elementary trimolecular reaction of the general form
\begin{equation}\label{eq:example_trimolecular_reaction}
    \ce{A + B + C -> C},
\end{equation}
which we suppose is mediated by a proximity-based reaction condition akin to those developed by Kearney and Flegg in \cite{our_first_paper}.\par

To describe the reaction condition, we follow \textcite{flegg2016smoluchowski} and generalise the notion of reactant proximity to trimolecular reactions by introducing two new coordinates, $\pos{\eta}{1}$ and $\pos{\eta}{2}$. We will discuss this coordinate system in greater detail in \cref{sec:math_background}, so for now it is sufficient to note that $\pos{\eta}{1}$ describes the separation between a particular pair of $A$ and $B$ molecules. While $\pos{\eta}{2}$ captures the relative proximity of a third molecule $C$ to this pair. Due to the radial symmetry of the system, proximity-based reaction conditions for the reaction in Equation \cref{eq:example_trimolecular_reaction} can be expressed naturally as functions of the radial coordinates $r_i \equiv ||\pos{\eta}{i}||$. Smoluchowski's original bimolecular condition can be extended in a straightforward manner to our trimolecular reaction by considering a reaction boundary that extends only a small distance in both $r_1$ and $r_2$. In this case, the system will reach a pseudo-equilibrium that returns trimolecular mass-action kinetics. However, a trimolecular system contains an additional spatial degree of freedom, namely $r_2$, when compared to a bimolecular system, which offers an opportunity to consider more elaborate reaction conditions. Specifically, we can consider reaction conditions where the separation, $r_1$, at which the pair of $A$ and $B$ molecules reacts is mediated by the distance, $r_2$, to the molecule of $C$ that is closest to this pair. In this case, the reaction condition still manifests as an absorbing boundary $\partial \Omega_R$ that is always small in the $\pos{\eta}{1}$ direction, meaning molecules of $A$ and $B$ must be close if they are to react, but now can be comparatively long in the $\pos{\eta}{2}$ direction so the molecule of $C$ involved in the reaction is not necessarily close to $A$ and $B$. \par

The expected reaction rate for the trimolecular system is calculated in the same way as it is for a bimolecular system and is determined by the total flux of the probability density across $\partial \Omega_R$. However, now we are concerned with the probability density for finding the state, where each state consists of the pair of $A$ and $B$ molecules and a particular molecule of $C$, associated with the $C$ molecule with the minimum $r_2$ value. That is, we must restrict our attention to the state that contains the molecule of $C$ that is closest to $A$ and $B$, since this state will be the first to satisfy the reaction condition. For all but the simplest choices of $\partial \Omega_R$, determining this probability density analytically is very difficult, and thus it is not possible to derive an analytic expression for the reaction rate in general. To avoid this issue, \textcite{our_first_paper} (implicitly) utilised singular perturbation theory to derive a leading-order solution for the steady-state probability density near $\partial \Omega_R$. They then demonstrated, via a particle-based simulation, that the total flux across this boundary can be approximated by just the diffusive flux of this leading-order solution in the $\pos{\eta}{1}$ direction. This approximation yields an analytic expression, analogous to Smoluchowski's in Equation \cref{eq:smol_steady_rate_constant}, that relates the reaction rate to the reaction boundary and thus allows us to match $\partial \Omega_R$ to any given reaction rate.\par

Unfortunately, this approximation will yield a reaction boundary that does not exactly reproduce the desired reaction rate when incorporated into a simulation, and instead small errors are observed between the simulated and desired rates. These errors are the result of neglecting second-order contributions to the diffusive flux in 
the $\pos{\eta}{1}$ direction across $\partial \Omega_R$ and entirely neglecting the flux across this boundary in the $\pos{\eta}{2}$ direction that arises from the relative motion of $C$ molecules. That is, the simulated rate reflects the exact flux across $\partial \Omega_R$ which in general differs from what we are led to expect when adopting the approximate reaction rate. In this article, we seek to address these shortcomings by examining the accuracy of this approximation for the functional form of $\partial \Omega_R$ considered by \cite{our_first_paper}. \par 

We begin in \cref{sec:math_background} by deriving a nonlinear partial integro-differential equation that governs the evolution of the probability density to find the state associated with the closest molecule of $C$ and whose analytic solution is unknown for arbitrary $\partial \Omega_R$. Although we consider its direct application to the optimisation of particle-based simulations, this equation has broader interest in applied mathematics. This work is, in part, a continuation of Redner and Ben-Avraham's \cite{SRedner_1990} in which they examined the distance of the closest molecule to a static trap. Redner and Ben-Avraham generalised the three-dimensional Hertz distribution \cite{Hertz-distribution}, which assumes that both the trap and surrounding molecules are static, but only considered a constant reaction radius. In contrast, we consider the case where this radius varies with $r_2$. Moreover, since all reactants undergo diffusion it is tempting to compare our results to the case where the trap is also mobile, which has been studied in lower dimensional systems \cite{mobile_trap, mobile_trap_later}, however, such studies are typically concerned with the correlations that arise between reactants \cite{scaling_approach_recombination_processes, spatial_structure_diffusion_limited, donev1999generalized}, whereas we restrict our attention to the mean-field behaviour of the system.\par 

We use singular perturbation theory \cite{KADALBAJOO2003371,neu2015singular} to derive a leading-order solution for probability density to find the closest molecule, and demonstrate how to use this solution to construct the aforementioned approximation to the steady-state reaction rate. Ideally, we would continue with this analysis to derive the higher-order corrections to the solution, and hence to the reaction rate, but the particular form of the problem makes this challenging. Thus, to quantify these corrections we turn to numerical methods, which are the other `principal approach for solving singular perturbation problems' \cite{KADALBAJOO2003371} and are commonly used to investigate singularly perturbed reaction diffusion systems \cite{Roos_robust_numerical_methods, doi:10.1137/110837784,KhariKartikay2022Aent,}. Since both the probability density and the flux of this quantity are of interest, in \cref{sec:numerical_methods} we adopt a mixed-primal variational formulation of the problem and present the corresponding finite element discretisation, which uses a continuous interior penalty scheme for the primal equation. The associated finite element solutions are then used to probe the higher-order corrections to the probability density and the reaction rate in \cref{sec:finite_element_solutions}. Finally, in \cref{sec:reactive_boundary_corrections} we demonstrate how to alter the reaction boundary so that it more accurately reproduces the desired reaction rate, enabling the development of more accurate particle-based simulations without sacrificing the convenience afforded by an analytic expression.

\section{The evolution of the closest molecule}\label{sec:math_background}
To model the elementary trimolecular reaction in Equation \cref{eq:example_trimolecular_reaction} we consider a domain $\Omega$ of volume $V$, where $V$ is finite but very large, that contains the three distinct chemical species, $A$, $B$ and $C$, whose molecules are all initially well mixed (distributed uniformly at random) within $\Omega$. It is convenient to assume that the system contains a single molecule of $A$ and $B$, and $N_C = cV$ molecules of $C$ where $c$ is a well-mixed concentration of $C$ molecules. In essence, we are assuming that distinct pairs of $A$ and $B$ molecules can be treated independently, which, as in the bimolecular case, is a good approximation so long as molecules only interact on small spatial scales. Furthermore, we assume that the molecules diffuse independently within the volume.
Under these assumptions, the system contains $N_C$ distinct states, where each state contains the single molecule of $A$, the single molecule of $B$ and a particular $C$ molecule.\par

To understand the dynamics of this system, it is instructive to initially restrict our attention to a single of the $N_C$ distinct states. We will use $D_i$ and $\mathbf{x}_i$, where $i = 0, 1, 2$, to denote the diffusion constant and the $3$-dimensional position of the molecules of $A$, $B$ and $C$ respectively. Since we are primarily concerned with the molecules' relative proximity it is convenient to transform the coordinate system into diffusive Jacobi coordinates or separation coordinates \cite{flegg2016smoluchowski},
\begin{subequations}
\begin{align}
    \pos{\eta}{0} &= \posbar{x}{2} \quad \text{and} \quad \label{eq:eta1_def}\\
    \pos{\eta}{i} &= \pos{x}{i} - \posbar{x}{i-1}, \quad i = 1,2,  \quad \text{where,}\label{eq:etai_def}\\
    \posbar{x}{i} &= \frac{\sum^i_{j=0} \pos{x}{j}D_j^{-1}}{\sum^i_{k=0}D_k^{-1}},\label{eq:x_bari_def}
\end{align}\end{subequations}
is the \textit{centre of diffusion} of the first $i+1$ molecules and is analogous to the centre of mass except that the positions are weighted by their inverse diffusion coefficients rather than their masses.\par

\begin{figure}
\centering
\includegraphics[width=0.5\linewidth]{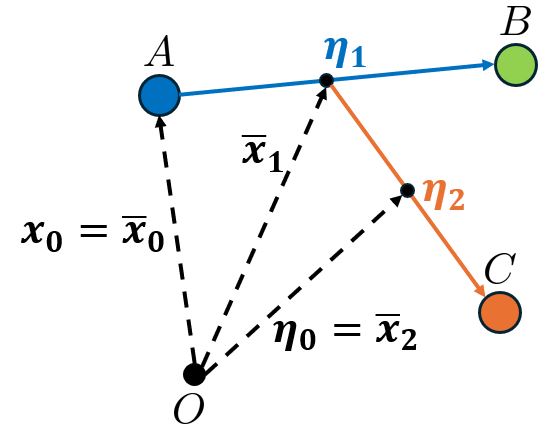}
\caption{The diffusive Jacobi coordinates for a system containing three molecules. The coordinate $\pos{\eta}{1}$ describes the separation between the first two molecules, labelled $A$ and $B$, while $\pos{\eta}{2}$ describes the separation of the third molecule $C$ from the centre of diffusion of $A$ and $B$ which we denote $\posbar{x}{1}$ and is defined in Equation \cref{eq:x_bari_def}. Finally, $\pos{\eta}{0}$ is the centre of diffusion of the three molecules ($\posbar{x}{2}$) and changes in this coordinate correspond to translations of the state.}
\label{fig:eta_coordinates}
\end{figure}

Since $\pos{\eta}{0}$ is the centre of diffusion of the three molecules, changes in it correspond to translations of the state within the domain. That is, varying $\pos{\eta}{0}$ does not change the relative proximity of the molecules in the current state, so we will only consider initial and boundary conditions that are independent of $\pos{\eta}{0}$. As a result, the joint probability density to find the unreacted state is also independent of $\pos{\eta}{0}$ and we will be concerned exclusively with $\pos{\eta}{1}$ and $\pos{\eta}{2}$. As shown in \cref{fig:eta_coordinates}, $\pos{\eta}{1}$ describes the separation between $A$ and $B$ as by definition $\posbar{x}{0} = \pos{x}{0}$, and $\pos{\eta}{2}$ describes the separation between $C$ and the centre of diffusion, $\posbar{x}{1}$, of $A$ and $B$. In addition, $\pos{\eta}{1}$ and $\pos{\eta}{2}$ can be shown to undergo independent linear diffusion with diffusion constants

\begin{align}\label{eq:eta_i_diff_coeff}
\hat{D}_1 &= D_{1} + \bar{D}_{0} \quad \text{and} \quad 
    \hat{D}_2 = D_{2} + \bar{D}_{1},                   
\end{align}
respectively, where $\bar{D}_j$ is the diffusion constant associated with $\posbar{x}{j}$
\begin{equation}\label{eq:x_bar_diff_coeff}
    \bar{D}_j = \frac{1}{\sum_{i=0}^{j} D_i^{-1}}.
\end{equation}\par

The joint probability density, $\prob{\pos{\eta}{}}{t}$, to find the molecules in an unreacted state at time $t$ is given by the diffusion equation
\begin{equation}\label{eq:general_prob_density_PDE}
    \frac{\partial \prob{\pos{\eta}{}}{t}}{\partial t} = \left[\sum_{i=1}^2 \hat{D}_i\hat{\nabla}^2_i \right]\prob{\pos{\eta}{}}{t},
\end{equation}
where $\boldsymbol{\eta} =  \left\{\pos{\eta}{1},\pos{\eta}{2}\right\}$ is the state vector for the system and $\hat{\nabla}^2_i$ denotes the Laplacian with respect to the coordinates of $\pos{\eta}{i}$. Since the system contains only a single molecule of $A$ and $B$, the entire system is absorbed when a reaction occurs. For the trimolecular reaction in Equation \cref{eq:example_trimolecular_reaction} a general proximity-based reaction condition can be defined as the absorbing boundary $\partial \Omega_R$ of the region
\begin{equation}\label{eq:generalised_inner_BC}
    \Omega_{R} = \left\{(\pos{\eta}{1}, \pos{\eta}{2} : r_1 \leq f(r_2)\right\},
\end{equation}
where $0 \leq f(r_2) \leq \sigmamax$ is a monotonically decreasing function of $r_2$ (recall that $r_i 
 = ||\pos{\eta}{i}||$) and $\sigmamax$ is a small positive constant relative to the diffusion coefficients $\hat{D}_i$ and the characteristic scale of the domain of $f$. In essence, the distance at which the molecules of $A$ and $B$ react is mediated by their relative proximity to the molecule of $C$ rather than being fixed as it is in the bimolecular case. Moreover, the monotonic decay of the reaction boundary reflects the physical notion that reactions to become more likely as the relative proximity between the three molecules reduces. In the case that the system contains a single molecule of $C$ its evolution is governed by Equation \cref{eq:general_prob_density_PDE} subject to the inner boundary condition
\begin{equation}\label{eq:general_prob_density_IB}
    P\left(\boldsymbol{\eta} \in \partial \Omega_R,t\right) = 0.
\end{equation}
However, in general this particular state won't be the first to cross $\partial \Omega_R$ and so to model the full system, we must consider $N_C$ distinct states; one for each combination of $A$, $B$ and a particular $C$ molecule.\par 

\begin{figure}[H]
    \centering
    \includegraphics[width=0.75\linewidth]{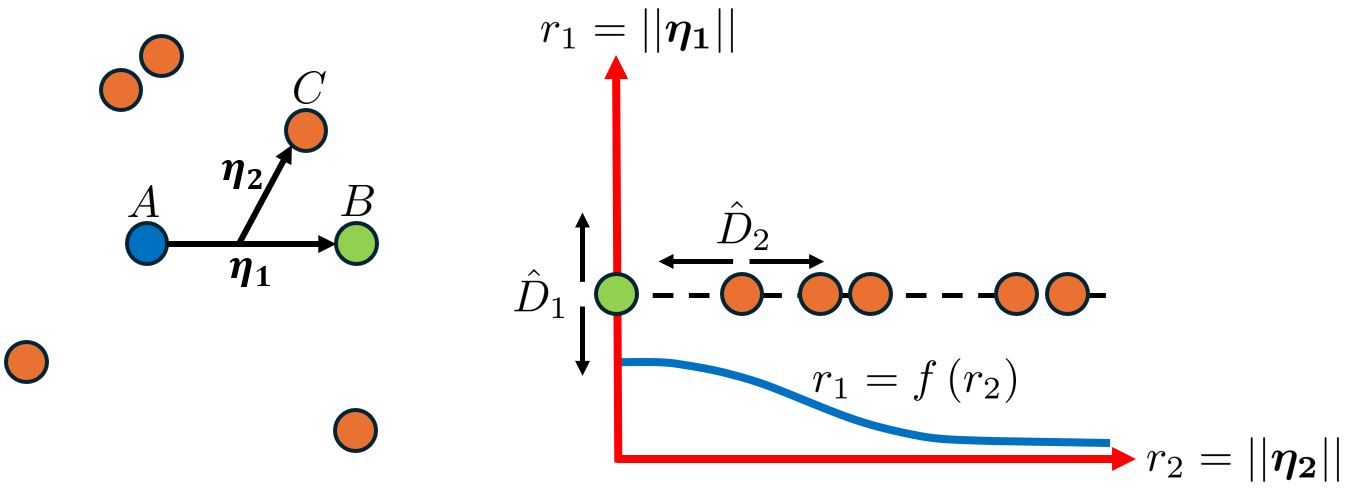}
    \caption{The state space for a system that contains a single molecule of $A$ (the blue point), a single molecule of $B$ (the green point) and five molecules of $C$ (the orange points). Each combination of $A$ 
 and $B$ with a particular molecule of $C$ gives a point within the state space on the right, shown by the orange points. The position of the $B$ molecule relative to the $A$ molecule defines a manifold of constant $\pos{\eta}{1}$ that all states lie on. This manifold diffuses in $\pos{\eta}{1}$ in accordance with Equation \cref{eq:diffusion_eta1} and the states on the manifold diffuse independently in $\pos{\eta}{2}$ according to Equation \cref{eq:diffusion_eta2}. The first state that crosses the inner boundary, depicted by the blue curve and defined in Equation \cref{eq:generalised_inner_BC}, will be absorbed.}
    \label{fig:state_diffusion}
\end{figure}

For any one particular state, $\pos{\eta}{1}$ and $\pos{\eta}{2}$ diffuse independently, so that the joint probability density can be written as, $P\parentheses{(}{\boldsymbol{\eta},t}{)} = P_1\parentheses{(}{\pos{\eta}{1},t}{)}P_2\parentheses{(}{\pos{\eta}{2},t}{)}$, where the probability densities $P_1$ and $P_2$ evolve according to
\begin{subequations}
\begin{align}
    \parentheses{[}{\frac{\partial}{\partial t} - \hat{D}_1\hat{\nabla}^2_1}{]}P_1\parentheses{(}{\pos{\eta}{1},t}{)} &\equiv\mathcal{L}_1P_1\parentheses{(}{\pos{\eta}{1},t}{)} = 0 \quad \text{and} \label{eq:diffusion_eta1}\\
    \parentheses{[}{\frac{\partial}{\partial t} - \hat{D}_2\hat{\nabla}^2_2}{]}P_2\parentheses{(}{\pos{\eta}{2},t}{)} &\equiv\mathcal{L}_2P_2\parentheses{(}{\pos{\eta}{2},t}{)} = 0\label{eq:diffusion_eta2}.
\end{align}\end{subequations}
Here $\mathcal{L}_1$ and $\mathcal{L}_2$ are diffusion operators on the $3$-dimensional spaces defined by $\pos{\eta}{1}$ and $\pos{\eta}{2}$ respectively. Since every state shares the same molecule of $B$, all states lie on manifolds of constant $\pos{\eta}{1}$ and the specific instance of $\pos{\eta}{1}$ diffuses according to Equation \cref{eq:diffusion_eta1} as shown in \cref{fig:state_diffusion}. On this manifold, states diffuse independently in the 3-dimensional space defined by $\pos{\eta}{2}$ in accordance with Equation \cref{eq:diffusion_eta2}. The first state incident on the inner boundary $\partial \Omega_R$ will cause the entire system to be absorbed, and since $f$ is monotonic, this state will be the one with the minimum $r_2$. That is, we wish to know the dynamics of the state that contains the molecule of $C$ that is closest to the centre of diffusion of $A$ and $B$.\par 

If the system is well mixed, we can assume that the $\pos{\eta}{2}$ coordinates of the $N_C$ states are uniformly and independently distributed throughout $\Omega$. Now consider a particular molecule of $C$ and let $\pos{Z}{2}$ denote the event that this particular molecule is associated with the minimum $r_2$ when compared with any other $C$ molecule in the system. Since the molecules are well mixed, the probability that an arbitrary $C$ molecule has the minimum $r_2$ at time $t$ is
\begin{equation}\label{eq:prob_eta1}
    \text{Pr}\left(\pos{Z}{2},t\right) = \frac{1}{N_C}.
\end{equation}
The probability that a particular molecule of $C$ has the minimum $r_2$ at time $t$ given it has a known $\pos{\eta}{2}$, is equal to the probability that all the other $N_C-1$ molecules of $C$ lie outside a sphere, $V_2$, of radius $r_2$ centred on the origin
\begin{equation}
\label{eq:Z2_given_eta1}
    \text{Pr}\left(\pos{Z}{2},t|\pos{\eta}{2}\right) = \left[1 - \int_{V_2}\probtwo[]{\pos{\eta'}{2}}{t}dV_2'\right]^{N_C-1},
\end{equation}
where $dV_2'$ is an elemental volume for the coordinates of $\pos{\eta'}{2}$. The probability density, $\probGiven[]{\pos{\eta}{2}}{t}{\pos{Z}{2}}$, for finding a particular molecule of $C$ at $\pos{\eta}{2}$ at time $t$, given that it has the smallest $r_2$ of any molecule of $C$ in the system, can be found by using Bayes' theorem in combination with Equations \cref{eq:prob_eta1} and \cref{eq:Z2_given_eta1},
\begin{equation}\label{eq:bayes_theorem_sub}
     \probGiven[]{\pos{\eta}{2}}{t}{\pos{Z}{2}} = N_C\probtwo[]{\pos{\eta}{2}}{t}\left[1-\int_{V_2}\probtwo[]{\pos{\eta'}{2}}{t}dV_2'\right]^{N_C-1}.
\end{equation}
\par

The system is very large and, in the limit that $V$ - and hence $N_C$ - tends to infinity, $\prob[]{\pos{Z}{2}}{t}$ goes to zero in accordance with Equation \cref{eq:prob_eta1}. To ensure we take the appropriate limit in Equation \cref{eq:bayes_theorem_sub} we define the scaled probability densities
\begin{subequations}
\begin{align}
    \PhiFunc &= \frac{\probGiven[]{\pos{\eta}{2}}{t}{\pos{Z}{2}}}{c} \quad \text{and} \label{eq:Phi_definition}\\
    \phiFunc &= \frac{N_C\probtwo[]{\pos{\eta}{2}}{t}}{c} \label{eq:phi_definition},
\end{align}\end{subequations}
which when substituted into Equation \cref{eq:bayes_theorem_sub} give,
\begin{equation}\label{eq:PhiFunc_integral}
    \PhiFunc = \phiFunc\left[1-\frac{c}{N_C}\int_{V_2}\phiFuncPrime dV_2'\right]^{N_C-1}.
\end{equation}
Taking the limit $N_C \rightarrow \infty$ in Equation \cref{eq:PhiFunc_integral} gives
\begin{equation}\label{eq:Phi_definition_in_limit}
    \PhiFunc = \phiFunc \text{exp}\left(-c\int_{V_2}\phiFuncPrime dV_2'\right),
\end{equation}
and then using Equation \cref{eq:diffusion_eta2}, we find that $\PhiFunc$ evolves according the advection-diffusion equation
\begin{equation}\label{eq:Phi_evolution}
    \frac{\partial \PhiFunc}{\partial t} = \hat{D}_2\hat{\nabla}^2_2\PhiFunc + \hat{D}_2\hat{\nabla}_2 \cdot \parentheses{(}{4\pi r_2^2 c\phi\left(\pos{\eta}{2}, t\right)\PhiFunc\pos{\hat{r}}{2}}{)},
\end{equation}
where $\pos{\hat{r}}{2}$ is the unit outward facing normal vector of a sphere of radius $r_2$ centred on the origin. A derivation of this result can be found in \cref{sec:AppendixA}.\par 

Equation \cref{eq:Phi_evolution} can be written entirely in terms of $\PhiFunc$ by noting that $\probGiven[]{\pos{Z}{2}}{t}{\pos{\eta}{2}}$ is also equal to 
the probability that the closest molecule is not inside a sphere, $V_2$, of radius $r_2$ centred on the origin. That is, in addition to Equation \cref{eq:Z2_given_eta1} we must also have
\begin{equation}
    \text{Pr}\left(\pos{Z}{2},t|\pos{\eta}{2}\right) = 1 - \int_{V_2}\probGiven[]{\pos{\eta'}{2}}{t}{\pos{Z}{2}}dV_2' = \int_{r_2}^\infty \probGiven[]{\pos{\eta'}{2}}{t}{\pos{Z}{2}}4\pi r_2'^2 \,\mathrm{d}r_2',
\end{equation}
Bayes Theorem then gives
\begin{equation}
    \probtwo[]{\pos{\eta}{2}}{t} = \frac{\probGiven[]{\pos{\eta}{2}}{t}{\pos{Z}{2}}\text{Pr}\left(\pos{Z}{2},t\right)}{\text{Pr}\left(\pos{Z}{2},t|\pos{\eta}{2}\right)},
\end{equation}
from which we conclude
\begin{equation}\label{eq:alternative_phi_expression}
    c\phiFunc = \frac{\PhiFunc}{\int_{r_2}^\infty \Phi\left(\pos{\eta'}{1},t\right) 4\pi r_2'^2 \,\mathrm{d}r_2'}.
\end{equation}
Substituting Equation \cref{eq:alternative_phi_expression} into Equation \cref{eq:Phi_evolution}, and then using Equation \cref{eq:Phi_definition}, we arrive at the governing equation for $\probGiven[]{\pos{\eta}{2}}{t}{\pos{Z}{2}}$
\begin{equation}\label{eq:P_closest_eta1_evolution}
    \frac{\partial \probGiven[]{\pos{\eta}{2}}{t}{\pos{Z}{2}}}{\partial t} = \hat{D}_2\hat{\nabla}^2_2\probGiven[]{\pos{\eta}{2}}{t}{\pos{Z}{2}} + \hat{D}_2\hat{\nabla}_2 \cdot \left(\frac{r_2^2 \left(\probGiven[]{\pos{\eta}{2}}{t}{\pos{Z}{2}}\right)^2\pos{\hat{r}}{2}}{\int_{r_2}^\infty \probGiven[]{\pos{\eta'}{2}}{t}{\pos{Z}{2}}r_2'^2 \,\mathrm{d}r_2'}\right).
\end{equation}
Finally, we note that since $\pos{\eta}{1}$ diffuses independently of $\pos{\eta}{2}$ the evolution of the joint probability density, $\probGiven[]{\boldsymbol{\eta}}{t}{\pos{Z}{2}}$, for finding the separation of the state with the minimum value of $r_2$ is governed by,
\begin{equation}\label{eq:governing_min_eta1_PDE}
    \frac{\partial \probGiven[]{\boldsymbol{\eta}}{t}{\pos{Z}{2}}}{\partial t} =  \left[\sum_{i=1}^2 \hat{D}_i\hat{\nabla}^2_i \right]\probGiven[]{\boldsymbol{\eta}}{t}{\pos{Z}{2}} + \hat{D}_2\hat{\nabla}_2 \cdot \left(\frac{r_2^2 \left(\probGiven[]{\boldsymbol{\eta}}{t}{\pos{Z}{2}}\right)^2\pos{\hat{r}}{2}}{\int_{r_2}^\infty P(\pos{\eta}{1}, \pos{\eta'}{2},t|\pos{Z}{2})r_2'^2 \,\mathrm{d}r_2'}\right).
\end{equation}
The linear diffusion terms in Equation \cref{eq:governing_min_eta1_PDE} describe the independent Brownian motion of $\pos{\eta}{1}$ and $\pos{\eta}{2}$ while the advection towards $\pos{\eta}{2}=\mathbf{0}$ represents the flux of the likelihood that the $C$ molecule with the second-smallest $||\pos{\eta}{2}||$ value diffuses over the sphere of radius $r_2$ set by the current closest $C$ molecule. In other words, it accounts for the fact that the molecule of $C$ that is the second-closest to the origin can diffuse inwards and become the closest molecule to the origin.\par

In the bimolecular case considered by Smoluchowski, the probability density relates to the likelihood of finding any molecule of $B$ at a particular position, while in the trimolecular case we have been careful to restrict our attention to the state associated with the closest molecule of $C$. This disparity arises because we have assumed that our trimolecular reaction boundary $\partial \Omega_R$ extends a very small distance, at most $\sigmamax$, in the $\pos{\eta}{1}$ direction, but may be comparatively long in the $\pos{\eta}{2}$ direction. In the extreme case where $f(r_2) = \sigmamax$ any molecule of $C$ will satisfy the reaction condition regardless of its proximity to $A$ and $B$ once $r_1 \leq \sigmamax$, and we must specify a particular molecule of $C$ to be involved in the reaction to avoid this ambiguity. Moreover, particle-based simulations of reaction-diffusion systems often update the position of molecules using finite time steps \cite{brownian_dynamics_with_hydro,brownianDynamics, lipkova2011analysis} which means that molecules are moved using small discrete displacements rather than continuously. If a pair of $A$ and $B$ molecules is such that their separation $r_1$ is just larger than $f(r_2)$ they can `jump' through the reaction boundary during the next time step, so that $r_1 < f(r_2)$ following the position updates. At this point, there can again be multiple molecules of $C$ close enough to satisfy the reaction condition, even if $f$ is not constant. In theory, the same problem exists when simulating bimolecular reactions and the distribution for the closest molecule to Smoluchowski's reaction boundary can be derived analytically \cite{SRedner_1990}. However, in practice $\sigma$ is sufficiently small that the likelihood that two molecules of $B$ are close enough to $A$ to react at any given moment is insignificant. Thus, unlike the trimolecular case, it is always clear which molecule of $B$ should react with a particular molecule of $A$.\par

To determine a unique solution for $\probGiven[]{\boldsymbol{\eta}}{t}{\pos{Z}{2}}$ we require appropriate initial and boundary conditions. To this end we assume that initially all the molecules are distributed uniformly at random throughout $\Omega$ so that the initial condition is given by the boundary-free steady-state solution of Equation \cref{eq:governing_min_eta1_PDE}
\begin{equation}\label{eq:governing_min_eta1_IC}
    \probGiven[]{\boldsymbol{\eta}}{0}{\pos{Z}{2}} = \frac{c}{V}\text{exp}\left(\frac{-4\pi c r_2^3}{3}\right).
\end{equation}
Moreover, we know that the system will be absorbed when the state associated with the closest molecule of $C$ crosses $\partial \Omega_R$, which leads to the inner boundary condition
\begin{equation}\label{eq:governing_min_eta1_IB}
    \probGiven[]{\boldsymbol{\eta} \in \partial \Omega_R}{t}{\pos{Z}{2}} = 0.
\end{equation}
Notice that this condition reduces to Equation \cref{eq:general_prob_density_IB} when the domain contains only a single molecule of $C$. Infinitely far from $\partial \Omega_R$ in both the $\pos{\eta}{1}$ and $\pos{\eta}{2}$ direction, the probability density is unperturbed by the absorption of states on $\partial \Omega_R$ so that
\begin{subequations}
\begin{align}\label{eq:governing_min_eta1_r0OBC}
    \lim_{r_1 \to \infty}\probGiven[]{\boldsymbol{\eta}}{t}{\pos{Z}{2}}  &= \frac{c}{V}\text{exp}\left(\frac{-4\pi c r_2^3}{3}\right) \quad \text{and}\\
    \lim_{r_2 \to \infty}\probGiven[]{\boldsymbol{\eta}}{t}{\pos{Z}{2}}  &= \lim_{r_2 \to \infty} \frac{c}{V}\text{exp}\left(\frac{-4\pi c r_2^3}{3}\right) = 0.\label{eq:governing_min_eta1_r1OBC}
\end{align}\end{subequations}
Finally, the flux of the probability density in the $\pos{\eta}{2}$ direction should vanish at $r_2 = 0$ which amounts to requiring that
\begin{equation}\label{eq:governing_min_eta1_r1IBC}
\frac{\partial \probGiven[]{\boldsymbol{\eta}}{t}{\pos{Z}{2}}}{\partial r_2} = 0 \quad \text{at } r_2 = 0.
\end{equation}
\par

Equation \cref{eq:governing_min_eta1_PDE}, subject to the conditions in Equations \cref{eq:governing_min_eta1_IC}--\cref{eq:governing_min_eta1_r1IBC}, describes the evolution of the state that will cause the system to be absorbed. However, solving this PDE when $\partial \Omega_R$ is chosen arbitrarily is very difficult. If $\partial \Omega_R$ extends a small distance in both $r_1$ and $r_2$ (rather than just in $r_1$ as assumed in Equation \cref{eq:generalised_inner_BC}) then the advection term in Equation \cref{eq:governing_min_eta1_PDE} becomes negligible, and the problem reduces to the trimolecular generalisation of Smoluchowski's reaction condition considered by \textcite{flegg2016smoluchowski}, for which an exact solution can be found. Outside of this regime, we are not aware of an exact solution, but it is possible to construct a leading-order solution using singular perturbation methods.\par 

Similar to the bimolecular case considered by Smoluchowski, the full time-dependent solution of Equation \cref{eq:governing_min_eta1_PDE} will rapidly converge to a steady state, since $\partial \Omega_R$ only extends a small distance in the $\pos{\eta}{1}$ direction \cite{comprehensiveChemicalKinetics}.
Thus, we will focus on computing the steady-state solution of the PDE, which will be radially symmetric in both $\pos{\eta}{1}$ and $\pos{\eta}{2}$ because the boundary conditions (Equations \cref{eq:governing_min_eta1_IC}--\cref{eq:governing_min_eta1_r1IBC}) only depend on the corresponding radial coordinates $r_1$ and $r_2$. That is, the steady-state solution of Equation \cref{eq:governing_min_eta1_PDE}, which we will denote $\steadyProb$ henceforth, will only depend on $r_1$ and $r_2$.\par

In the limit that $\sigmamax \to 0$, $\partial \Omega_R$ becomes a sphere of radius $0$ and, in effect, the absorbing boundary vanishes. In this case, the steady-state solution to Equation \cref{eq:governing_min_eta1_PDE} is the initial condition from Equation \cref{eq:governing_min_eta1_IC}. However, this solution has no $r_1$ dependence and so will never satisfy the boundary condition on $\partial \Omega_R$ for any finite value of $\sigmamax$. This behaviour indicates that there exists a boundary layer \cite{Schlichting2017} near $\partial \Omega_R$. That is, in the vicinity of the inner boundary, we expect the steady-state solution to change rapidly from the boundary-free steady-state solution to one that satisfies the boundary condition.\par 

To determine the solution near $\partial \Omega_R$, we seek a radial steady-state solution, $\probbar{u_1}{u_2} = \steadyProb$, where
\begin{equation}
    u_1 = \frac{r_1}{\sigmamax} \quad \text{and} \quad u_2 = r_2,
\end{equation}
are scaled radial coordinates. With this change of variables, the steady-state form of Equation \cref{eq:governing_min_eta1_PDE}
becomes
\begin{equation}\label{eq:perturbation_pde}
     \frac{\hat{D}_1}{u_1^2}\frac{\partial}{\partial u_1}\left(u_1^2\frac{\partial \probbar{u_1}{u_2}}{\partial u_1}\right) + \frac{\sigmamax^2\hat{D}_2}{u_2^2}\frac{\partial}{\partial u_2}\left(u_2^2\frac{\partial \probbar{u_1}{u_2}}{\partial u_2} + \frac{u_2^4 \left(\probbar{u_1}{u_2}\right)^2}{\int_{u_2}^\infty \probbar{u_1}{u_2'}u_2'^2 du_2'}\right) = 0,
\end{equation}
while the associated boundary conditions are:
\begin{subequations}
\begin{align}
    \probbar{u_1}{u_2} &= 0 \quad \text{on} \quad u_1 = F(u_2), \quad\text{where } F(u_2) = \frac{f(u_2)}{\sigmamax}\label{eq:perturbation_IBC},\\
    \lim_{u_1 \to \infty}\probbar{u_1}{u_2}  &= \frac{c}{V}\text{exp}\left(\frac{-4\pi c u_2^3}{3}\right),\label{eq:perturbation_u0OBC}\\
    \lim_{u_2 \to \infty}\probbar{u_1}{u_2} &= \lim_{u_2 \to \infty}\frac{c}{V}\text{exp}\left(\frac{-4\pi c u_2^3}{3}\right) = 0 \quad \text{and}
 \label{eq:pertrubation_u1OBC}\\
\frac{\partial \probbar{u_1}{u_2}}{\partial u_2} &= 0,   \quad \text{when }u_2=0.\label{eq:perturbation_u1IBC}
\end{align}\end{subequations}
Following the change of variables the singular nature of the PDE becomes plain. Very close to $\partial \Omega_R$ diffusion in the $\pos{\hat{r}}{1}$ direction dominates and there is comparatively little diffusive or advective motion in the $\pos{\hat{r}}{2}$ direction, where $\pos{\hat{r}}{i}$ is the unit radial basis vector for $\pos{\eta}{i}$. In the limit that $\sigmamax \to 0$ motion in the $\pos{\hat{r}}{2}$ direction can be ignored, and it is only 
necessary to consider diffusive flux in the $\pos{\hat{r}}{1}$ direction.\par

More formally, we assume the solution to Equation \cref{eq:perturbation_pde} can be expanded in a perturbation series of the form
\begin{equation}
    \probbar{u_1}{u_2} = \bar{P}_0\left(u_1,u_2\right) + \sigmamax^2\bar{P}_1\left(u_1,u_2\right) + O(\sigmamax^4).
\end{equation}
To derive the leading-order solution we substitute this perturbation series into Equation \cref{eq:perturbation_pde} and match the $O(1)$ terms which leads to
\begin{equation}\label{eq:P0_partial_solution}
    \bar{P}_0\left(u_1,u_2\right) = h(u_2)\left(1-\frac{F(u_2)}{u_1}\right),
\end{equation}
where $h$ is an arbitrary function of $u_2$. To determine $h$ we use the matching condition
\begin{equation}
    \lim_{u_1 \to \infty} \bar{P}_0\left(u_1,u_2\right) = \frac{c}{V}\exp\left(\frac{-4\pi c u_2^3}{3}\right),
\end{equation}
which amounts to requiring that in the bulk, far from $\partial \Omega_R$, the absorption of states on $\partial \Omega_R$ causes a negligible perturbation in the solution. Finally, reverting to our original coordinates $r_1$ and $r_2$ we recover the approximate steady-state solution
\begin{equation}\label{eq:P0_solution}
    P_0\left(r_1, r_2\right) = \frac{c}{V}\exp\left(\frac{-4\pi c r_2^3}{3}\right)\left(1-\frac{f(r_2)}{r_1}\right),
\end{equation}
which we will refer to as the leading-order solution henceforth.
\par

The steady-state reaction rate is the total flux of $\steadyProb$ over $\partial \Omega_R$. However, the outward facing unit normal vector to this boundary is 
\begin{equation}\label{eq:normal_vector}
    \pos{\hat{n}}{\text{out}} = -\left(1 + \left(\sigmamax\frac{d F(r_2)}{\,\mathrm{d}r_2}\right)^2\right)^{-\frac{1}{2}}\left(\pos{\hat{r}}{1} - \sigmamax\frac{d F(r_2)}{\,\mathrm{d}r_2}\pos{\hat{r}}{2}\right),
\end{equation}
hence, to first order in $\sigmamax$, only the diffusive flux in the $\pos{\hat{r}}{1}$ direction contributes.
That is, the steady-state reaction rate can be approximated by
\begin{equation}\label{eq:closest_reaction_rate}
    K_1(c) = \int_0^\infty \left(\hat{D}_1\nabla_1 P\left(\boldsymbol{\eta}|\boldsymbol{Z}_2\right)\cdot 4\pi (f(r_2))^2 \hat{\boldsymbol{r}}_1\right) 4\pi r_2^2 \,\mathrm{d}r_2,
\end{equation}
which can be simplified to
\begin{equation}\label{eq:closest_rate_old}
    K_1(c) = \int_0^\infty \frac{4\pi \hat{D}_1 f(r_2)}{V}c\text{exp}\left(\frac{-4\pi cr_2^3}{3}\right)4\pi r_2^2 \,\mathrm{d}r_2,
\end{equation}
using Equation \cref{eq:P0_solution} and neglecting the $O(\sigmamax^2)$ terms.\par

Equation \cref{eq:closest_rate_old} matches the expression used by \textcite{our_first_paper} to determine the function $f(r_2)$, which defines the boundary $\partial \Omega_R$ required to reproduce a given reaction rate $K_1(c)$. However, to arrive at this expression, several contributions to the flux across $\partial \Omega_R$ must be ignored. Firstly, there is an $O(\sigmamax^2)$ correction to the leading-order solution and hence to the diffusive flux in the $\pos{\hat{r}}{1}$ direction. In addition, there is an $O(\sigmamax^2)$ contribution to the flux that arises from the diffusive and advective fluxes in the $\pos{\hat{r}}{2}$ direction. Neglecting these contributions means that the functional form derived for $\partial \Omega_R$ will not be exact, resulting in small errors in the reaction rate reproduced by particle-based simulations that use this reaction boundary. Although these errors were known to Kearney and Flegg in \cite{our_first_paper}, they did not examine them in detail because it was difficult to isolate them from other sources of error present in their simulations. Ideally, we would quantify these corrections analytically by continuing our singular perturbation analysis of Equation \cref{eq:perturbation_pde}. Unfortunately, retaining the $O(\sigmamax^2)$ terms in this equation leaves us unable to find an analytic expression that satisfies the associated boundary conditions. Therefore, to quantify the expected $O(\sigmamax^2)$ corrections, we seek a finite element solution to the problem that can be compared to the leading-order solution in Equation \cref{eq:P0_solution}.

\section{Mixed-primal formulation and finite element discretisation}\label{sec:numerical_methods}
The quantity of primary concern for particle-based simulations is the steady-state reaction rate, which corresponds to the flux of $\steadyProb$ over the reaction boundary $\partial \Omega_R$. One approach to compute this flux is to use a reconstruction from the finite element solution to the original form of the PDE (see, e.g., \cite{local_flux_reconstruction}), but we opt for using a mixed-primal reformulation of the problem (which might be more appropriate to preserve the divergence structure \cite{gatica2014simple}), and solve for the flux
\begin{equation}\label{eq:s_definition}
    \s = -\left(\nabla \steadyProb + \frac{r_2^2 \steadyProb^2}{\Q}\hat{\mathbf{r}}_2\right),
\end{equation}
directly, where 
\begin{equation}
\label{eq:radial-grad}     \nabla \steadyProb = \frac{\partial \steadyProb}{\partial r_1}\hat{\mathbf{r}}_1 + \frac{\partial \steadyProb}{\partial r_2}\hat{\mathbf{r}}_2.\end{equation}
Notice that we have introduced the \textit{density moment}
\begin{equation}
    \Q = \int^\infty_{r_2} P(r_1, r_2')r_2'^2 \,\mathrm{d}r_2',  
\end{equation}
which enables us to avoid directly calculating the integral in the Equation \cref{eq:governing_min_eta1_PDE}. The steady-state is now governed by a coupled system of PDEs 
\begin{subequations}\label{eq:numerical_system}
\begin{align}
    \s + \left(\nabla \steadyProb + \frac{r_2^2 \steadyProb^2}{\Q}\hat{\mathbf{r}}_2\right) &= 0,\\  \nabla \cdot (D\s) &= 0, 
     \quad \text{and}\\
     \frac{\partial \Q}{\partial r_2} + r_2^2\steadyProb &= 0,
\end{align}\end{subequations}
where $D$ is a matrix of the diffusion coefficients such that,
\begin{equation}
    D = \begin{bmatrix}
\hat{D}_1 & 0\\
0 & \hat{D}_2
\end{bmatrix},
\end{equation}
and $\nabla \cdot \s$ denotes the radial divergence of $\s$, that is,
\begin{equation}\label{eq:radial-div}
    \nabla \cdot \s = \frac{1}{r_1^2}\frac{\partial}{\partial  r_1}\left(r_1^2\s \cdot \pos{\hat{r}}{1}\right) + \frac{1}{r_2^2}\frac{\partial}{\partial  r_2}\left(r_2^2\s \cdot \pos{\hat{r}}{2}\right).
\end{equation}
The corresponding steady-state reaction rate can be computed straightforwardly by evaluating
\begin{equation}\label{eq:total_flux_s}
   K(c) = \int_0^\infty D \s \cdot \pos{\hat{n}}{\text{out}}(4\pi f(r_2)r_2)^2 \,\mathrm{d}r_2, 
\end{equation}
where we recall that $\pos{\hat{n}}{\text{out}}$ is the outward facing unit normal vector to $\partial \Omega_R$.\par

In principle the boundary conditions for $\steadyProb$ remain unchanged except that the no-flux condition defined in Equation \cref{eq:governing_min_eta1_r1IBC} now defines a boundary condition on $\s$ 
\begin{equation}
    \s \cdot \pos{\hat{r}}{2} = 0\quad \text{at } r_2 = 0.
\end{equation}
However, we have assumed that the domain $\Omega$ is very large (effectively infinite) in defining these conditions, and representing such a large domain inherently limits the resolution of the finite element mesh. This is problematic because we must be able to resolve the rapid change in the solution expected within the boundary layer near $\partial \Omega_R$. Therefore, it is necessary to consider a domain amenable for the finite element discretisation 
\begin{equation}
    \Omega_{\text{FE}} = \left\{(r_1, r_2) : f(r_2) \leq r_1 \leq r_1^{\text{max}},0 \leq r_2 \leq r_2^{\text{max}}\right\},
\end{equation}
that is not necessarily sufficiently large for these boundary conditions to be valid. The inner boundary condition from Equation \cref{eq:governing_min_eta1_IC} transfers trivially to $\Omega_{\text{FE}}$,
\begin{equation}
    \steadyProb = 0 \quad \text{on } r_1 = f(r_2),
\end{equation}
but the remaining conditions in Equations \cref{eq:governing_min_eta1_r0OBC} and \cref{eq:governing_min_eta1_r1OBC} require more careful consideration.\par 

We expect the leading-order contribution to $\steadyProb$ to decay exponentially as $r_2$ increases. The rate of decay is expected to match that of the value on the boundary, given in Equation \cref{eq:governing_min_eta1_r1OBC}, and so the value imposed at $r_2 = r_2^{\text{max}}$ very quickly approaches zero. Therefore, the error that arises from imposing
\begin{equation}\label{eq:r2_max_condition}
    \steadyProb = \frac{c}{V}\text{exp}\left(\frac{-4\pi c r_2^3}{3}\right) \quad \text{at } r_2 = r_2^{\text{max}},
\end{equation}
for finite $r_2^{\text{max}}$, very quickly approaches zero as $r_2^{\text{max}}$ is increased. Adopting this condition implies that we must also have
\begin{equation}\label{eq:r2_max_Q_condition}
    \Q = \frac{1}{4\pi V}\text{exp}\left(\frac{-4\pi c r_2^3}{3}\right) \quad \text{at } r_2 = r_2^{\text{max}}.
\end{equation}\par

In contrast, as $r_1$ increases, we expect the leading-order solution to decay to the boundary condition in Equation \cref{eq:governing_min_eta1_r0OBC} with the reciprocal of $r_1$. That is, we expect the solution to converge to the outer solution much more slowly in the $\pos{\hat{r}}{1}$ direction than in the $\pos{\hat{r}}{2}$ direction and require $r_1^{\text{max}}$ to be significantly larger than $r_2^{\text{max}}$ for the boundary condition to be correct. Even if $r_1^{\text{max}}$ is $O(1/\sigma_\text{max})$ adopting the boundary condition in Equation \cref{eq:governing_min_eta1_r0OBC} incurs an $O(\sigma_\text{max}^2)$ error that arises from the $O(\sigmamax)$ contribution to $\steadyProb$. To avoid introducing this unnecessary error we use the leading-order solution to modify the boundary condition for finite $r_1^{\text{max}}$ and instead impose,
\begin{equation}\label{eq:r1_max_condition}
    \steadyProb = \frac{c}{V}\text{exp}\left(\frac{-4\pi c r_2^3}{3}\right)\left(1 - \frac{f(r_2)}{r_1}\right) \quad \text{at } r_1 = r_1^{\text{max}}.
\end{equation}
This boundary condition neglects any $O(\sigma_\text{max}^2)$ corrections to the solution; however, we expect these corrections to decay polynomially with increasing $r_1$ similar to the $O(\sigma_\text{max})$ term in the leading-order solution. Therefore, the error introduced by adopting this condition can be made negligible by making $r_1^{\text{max}}$ sufficiently large. \par


Before stating the finite element formulation, let us recall, for a sufficiently smooth vector-valued function $\btau=(\tau_1,\tau_2)^T$ and a scalar function $v$, the following integration by parts formula for radial coordinates. This employs the two-dimensional gradient from Equation \cref{eq:radial-grad}, usual two-dimensional divergence $\nabla_c\cdot \btau = \partial_{r_1} \tau_1 + \partial_{r_2} \tau_2$, radial divergence from Equation \cref{eq:radial-div} (noting that $\nabla \cdot \btau = \nabla_c\cdot \btau + \frac{2}{r_1} \tau_1 +  \frac{2}{r_2} \tau_2$), and usual integration by parts in  two-dimensions:
\begin{align}\label{eq:int-by-parts}
\nonumber\int_{\Omega_{\mathrm{FE}}} ( \nabla\cdot \btau) (v \omega)\, \mathrm{d}r_2\,\mathrm{d}r_1  & = \int_{\Omega_{\mathrm{FE}}} (\nabla_c\cdot \btau) (v\omega)\, \mathrm{d}r_2\,\mathrm{d}r_1 +  \int_{\Omega_{\mathrm{FE}}} 2\biggl(\frac{\omega }{r_1} \tau_1 +  \frac{\omega}{r_2} \tau_2\biggr) v\, \mathrm{d}r_2\,\mathrm{d}r_1 \\
\nonumber & =  \int_{\partial \Omega_{\mathrm{FE}}} (\btau\cdot \boldsymbol{n})(v \omega)\, \mathrm{d}s   - \int_{\Omega_{\mathrm{FE}}} (\btau\cdot \nabla v) \, \omega\, \mathrm{d}r_2\,\mathrm{d}r_1 - \int_{\Omega_{\mathrm{FE}}} (\btau\cdot \nabla \omega) v\, \mathrm{d}r_2\,\mathrm{d}r_1  \\
\nonumber & \qquad +  \int_{\Omega_{\mathrm{FE}}} 2\biggl(\frac{\omega }{r_1} \tau_1 +  \frac{\omega}{r_2} \tau_2\biggr) v\, \mathrm{d}r_2\,\mathrm{d}r_1 \\
& = \int_{\partial \Omega_{\mathrm{FE}}} (\btau\cdot \boldsymbol{n})(v \omega)\, \mathrm{d}s -  \int_{\Omega_{\mathrm{FE}}} (\btau\cdot \nabla v) \, \omega\, \mathrm{d}r_2\,\mathrm{d}r_1.
\end{align}
Here $\omega=(4\pi r_1r_2)^2$ denotes the weight associated with the radial coordinates $r_1$ and $r_2$, while $\partial \Omega_{\mathrm{FE}}$ denotes the boundary of $\Omega_{\mathrm{FE}}$ which we divide into four segments: $\Gamma_{\mathrm{left}}$, $\Gamma_{\mathrm{bot}}$, $\Gamma_{\mathrm{right}}$ and $\Gamma_{\mathrm{top}}$, defined by the regions $r_2 = 0$, $r_1 = f(r_2)$, $r_2 = r_2^{\text{max}}$ and $r_1 = r_1^{\text{max}}$ respectively.

Consider now a finite element mesh $\mathcal{T}_h$---an unstructured partition of the domain $\Omega_{\text{FE}}$ into triangles $K$---of maximum diameter $h=\{\max h_K: K \in \mathcal{T}_h\}$. Let us denote by $\mathcal{E}_h$ the set of all interior edges $e$ with length $h_e$ (and $\mathcal{E}_h^\partial$ denotes the boundary edges). In addition, we use $[\![ s]\!]_e$ to denote the jump of a scalar field $s$ across the edge $e$. Proceeding from a weak formulation of \eqref{eq:numerical_system} (derived from the integration by parts in radial coordinates \eqref{eq:int-by-parts} and using the mixed type boundary conditions -- for the flux essentially, and for the density and density moment naturally) and considering a mixed-primal Galerkin discretisation, the finite element formulation consists in finding the approximate unknowns of flux $\bs_h$, density $P_h$, and density moment $Q_h$ in the discrete spaces $\Sigma_h\times \Phi_h \times \Psi_h$, such that
\begin{subequations}\label{eq:fem}
    \begin{align}
      \int_{\Omega_{\mathrm{FE}}}\!\! \! \bs_h\cdot D\btau_h \omega\, \mathrm{d}r_2\,\mathrm{d}r_1  
       -  \int_{\Omega_{\mathrm{FE}}}\!\! \!\! \!P_h \nabla\cdot (D\btau_h) \omega\, \mathrm{d}r_2\,\mathrm{d}r_1  
      +   \int_{\Omega_{\mathrm{FE}}}\!\!  \!\biggl(\frac{r_2^2P_h^2}{Q_h}\pos{\hat{r}}{2}\biggr)\cdot (D\btau_h) \omega\, \mathrm{d}r_2\,\mathrm{d}r_1  & = - \int_{\Gamma_{\mathrm{right}}}\!\! \! \!\! \!\!\!P_{\mathrm{right}}D\btau_h\cdot\boldsymbol{n} \omega \mathrm{d}s \nonumber \\
      & \quad \, - \int_{\Gamma_{\mathrm{bot}}}\!\! \!\! \!P_{\mathrm{bot}}D\btau_h\cdot\boldsymbol{n} \omega \mathrm{d}s \nonumber \\
      & \quad \, - \int_{\Gamma_{\mathrm{top}}}\!\! \!\! \!P_{\mathrm{top}}D\btau_h\cdot\boldsymbol{n} \omega \mathrm{d}s \nonumber \\
      & \qquad \ \ \, \forall \btau_h \in \Sigma_h^0,\\
      -\int_{\Omega_{\mathrm{FE}}} v_h \nabla\cdot (D\bs_h) \omega\, \mathrm{d}r_2\,\mathrm{d}r_1 & = 0 \quad \forall v_h \in \Phi_h,\\
    \int_{\Omega_{\mathrm{FE}}} (r_2^2P_h- \nabla \cdot \pos{\hat{r}}{2} Q_h) w_h \omega\, \mathrm{d}r_2\,\mathrm{d}r_1 - 
      \int_{\Omega_{\mathrm{FE}}} (\nabla w_h \cdot \pos{\hat{r}}{2}) Q_h \omega\, \mathrm{d}r_2\,\mathrm{d}r_1 + \int_{\Gamma_{\mathrm{left}}} \!\!\!\!\!\pos{\hat{r}}{2} \cdot\boldsymbol{n} Q_hw_h\omega \mathrm{d}s & \nonumber \\
      + \sigma_{\mathrm{stab}} \sum_{e\in \mathcal{E}_h} \frac{h_e^2}{(k+1)^\alpha}\beta \int_e [\![ \nabla Q_h \cdot \boldsymbol{n}_e]\!]_e  [\![ \nabla w_h \cdot \boldsymbol{n}_e]\!]_e \omega \mathrm{d}s \nonumber \\
      + \sigma_{\mathrm{stab}}'\sum_{h_e\in \mathcal{E}_h^\partial} \bigl(h_e[\pos{\hat{r}}{2}\cdot\boldsymbol{n} ]^2_- + \frac{8}{h_e}\bigr) \int_eQ_hw_h \omega \mathrm{d}s & = \nonumber \\
      \int_{\Gamma_{\mathrm{right}}} \!\!\!\!\!\!\!\!\pos{\hat{r}}{2} \cdot\boldsymbol{n} Q_{\mathrm{right}} w_h\omega \mathrm{d}s + \sigma_{\mathrm{stab}}'\sum_{h_e\in \mathcal{E}_h^\partial} \bigl(h_e[\pos{\hat{r}}{2}\cdot\boldsymbol{n} ]^2_- + \frac{8}{h_e}\bigr)\int_e Q_{\partial}w_h \omega \mathrm{d}s   
      &  \qquad \ \ \, \forall w_h \in \Psi_h. \label{eq:Qh}
    \end{align}
\end{subequations}
The discrete trial and test spaces for the flux-density pair are, for given $k\geq 0$, Raviart--Thomas elements of degree $k$ and overall discontinuous and piecewise polynomials of degree $k$, while for the density moment we take overall continuous and piecewise polynomials of degree $k+1$ (see their precise definition for the Cartesian case in, e.g., \cite{ernguermond} and the modification for axisymmetric divergence in \cite{ervin2013approximation}, see also \cite{baird2021second,neilan2019low}): 
\begin{subequations}
    \begin{align}
  \Sigma_h &:= \{ \btau_h \in \mathbf{H}(\mathrm{div}_r;\Omega_{\mathrm{FE}}): \btau_h|_K \in \mathbb{RT}_k(K) \ \forall K\in\mathcal{T}_h, \ \btau_h\cdot \boldsymbol{n}|_{\Gamma_{\mathrm{left}}} = s_{\mathrm{left}} \},\\
  \Sigma_h^0 &:= \{ \btau_h \in \mathbf{H}(\mathrm{div}_r;\Omega_{\mathrm{FE}}): \btau_h|_K \in \mathbb{RT}_k(K) \ \forall K\in\mathcal{T}_h, \ \btau_h\cdot\boldsymbol{n}|_{\Gamma_{\mathrm{left}}} = 0 \},\\
 \Phi_h &:= \{v_h \in L^2_\omega(\Omega_{\mathrm{FE}}): v_h|_K \in \mathbb{P}_{k}(K) \ \forall K\in \mathcal{T}_h \},\\
  \Psi_h &:= \{w_h \in C^0(\overline{\Omega_{\mathrm{FE}}}): w_h|_K \in \mathbb{P}_{k+1}(K) \ \forall K\in \mathcal{T}_h\}.
    \end{align}
\end{subequations}
Here by $\mathbf{H}(\mathrm{div}_r;\Omega_{\mathrm{FE}})$ we mean all vector fields in $\mathbf{L}^2_\omega(\Omega_{\mathrm{FE}})$ such that their \emph{radial} divergence \eqref{eq:radial-div} is in the scalar version of the space, $L^2_{\omega}(\Omega_{\mathrm{FE}})$. In turn, for all $p>1$, by $L^p_{\omega}(\Omega_{\mathrm{FE}})$ we denote the weighted Lebesgue space of all measurable functions $v$ on $\Omega_{\mathrm{FE}}$ such that (see \cite{bernardi1999spectral}) 
\[ \| v\|^p_{L^p_\omega(\Omega_{\mathrm{FE}})} := \int_{\Omega_{\mathrm{FE}}} |v|^p \omega\, \mathrm{d}r_2\,\mathrm{d}r_1 < \infty.\]
Observe that in the discrete weak formulation we have tested the constitutive equation for the flux against a weighted test function $D\btau_h$ in order to obtain a block-symmetric structure in the weighted divergence operator. Note also, that since the equation for $Q_h$ is hyperbolic in nature, we employ a continuous interior penalty (C0IP) approach adopted from \cite{burman2007continuous}, that restores the otherwise suboptimal convergence expected when using the space $\Psi_h$. This method includes an additional integration by parts (that we conduct using the modified form \eqref{eq:int-by-parts}) and the interior jump penalisation associated with the normal gradients. Other possible choices are classical streamline upwind Petrov--Galerkin or streamline diffusion variants. More details can be found in, e.g., \cite{bochev2004stability,burman2007continuous}.\par 

We confirm experimentally the convergence of \eqref{eq:fem} to the unique weak solution of \eqref{eq:numerical_system} by using smooth manufactured solutions $P(r_1,r_2) = \sin(\pi r_2)\sin(\pi r_1)$, $Q(r_1, r_2) = \frac32 + \cos(\pi r_1r_2)$, and considering a manufactured right-hand side for the second and third equations of \eqref{eq:numerical_system}. We impose natural boundary conditions for $P_h$ on the right, bottom, and top parts of the boundary, essential boundary conditions for the flux on $\Gamma_{\mathrm{left}}$ and natural boundary conditions for $Q_h$ on the boundary $\Gamma_{\mathrm{right}}$ (the prescribed boundary values and source terms are not necessarily homogeneous but given by the manufactured exact density and density moment). We impose the latter conditions naturally since the upstream mechanism does not permit us to set up a boundary condition essentially on the right boundary (as this can be considered an outflow segment on which $\pos{\hat{r}}{2}\cdot \boldsymbol{n}>0$). This explains the additional terms on $\mathcal{E}_h^\partial$ suggested in \cite{burman2019primal}. A sequence of successively refined meshes is constructed, and on each level the approximate solutions are compared to the exact ones using the natural norms
\[e(\bs):=\|\bs-\bs_h\|_{\mathbf{H}(\mathrm{div}_r;\Omega_{\mathrm{FE}})}, \quad 
e(P):=\|P-P_h\|_{L^2_\omega(\Omega_{\mathrm{FE}})}, \quad 
e(Q):=\|Q-Q_h\|_{H_\omega(\Omega_{\mathrm{FE}})},
\]
and experimental convergence rates are computed as 
$r_{l\,+\,1}(\cdot) = \log(e_{l\,+\,1}(\cdot)/ e_{l}(\cdot))[\log(h_{l\,+\,1}/h_l)]^{-1}$, where $e_l$ denotes the error incurred at the refinement level with a mesh of size $h_l$. Here 
$$\|\bullet\|^2_{H_\omega(\Omega_{\mathrm{FE}})} = \|\bullet\|^2_{L^2_\omega(\Omega_{\mathrm{FE}})} + \int_{\Omega_{\mathrm{FE}}} |\pos{\hat{r}}{2}\cdot \nabla (\bullet)|^2  \omega\, \mathrm{d}r_2\,\mathrm{d}r_1 $$
is the norm in the graph space. 
The constants are taken as $\hat{D}_1 = 2, \hat{D}_2=1.5$, and we test the schemes with polynomial degrees $k=0$ and $k=1$. The numerical parameters associated with the C0IP discretisation \eqref{eq:Qh} are taken as $\alpha = \frac72$, $\beta = 1$, $\sigma_{\mathrm{stab}}=0.1$, and $\sigma_{\mathrm{stab}}'=k+1$.

\begin{table}[t!]
    \centering
    {\small\begin{tabular}{|rc|cccccc|}
    \hline
     DoF  &    $h$   &  $e(\bs)$   &  $r(\bs)$  &  $e(P)$   &  $r(P)$   & $e(Q)$ & $r(Q)$ \\  
     \hline 
     \multicolumn{8}{|c|}{Lowest-order method with $k=0$}\\
    \hline 
      33 & 0.7071 & 2.93e+01 &  $\star$ & 1.03e+00 &  $\star$ & 4.92e+00 &  $\star$ \\ 
   113 & 0.3536 & 1.73e+01 & 0.763 & 5.56e-01 & 0.888 & 7.08e+00 & -0.526 \\
   417 & 0.1768 & 9.96e+00 & 0.793 & 2.84e-01 & 0.969 & 1.84e+00 & 1.943 \\
  1601 & 0.0884 & 5.93e+00 & 0.747 & 1.43e-01 & 0.992 & 7.10e-01 & 1.375 \\
  6273 & 0.0442 & 3.72e+00 & 0.673 & 7.16e-02 & 0.997 & 3.35e-01 & 1.085 \\
 24833 & 0.0221 & 2.45e+00 & 0.604 & 3.59e-02 & 0.995 & 1.63e-01 & 1.036 \\
 98817 & 0.0110 & 1.67e+00 & 0.551 & 1.82e-02 & 0.981 & 8.13e-02 & 1.004 \\
\hline 
     \multicolumn{8}{|c|}{Scheme with $k=1$}\\
    \hline 
     97 & 0.7071 & 1.25e+01 & $\star$ & 1.92e-01 & $\star$ & 3.60e+00 & $\star$ \\
   353 & 0.3536 & 4.35e+00 & 1.522 & 7.00e-02 & 1.456 & 1.10e+00 & 1.717 \\
  1345 & 0.1768 & 1.41e+00 & 1.625 & 1.86e-02 & 1.916 & 1.81e-01 & 2.602 \\
  5249 & 0.0884 & 4.79e-01 & 1.557 & 4.70e-03 & 1.981 & 4.79e-02 & 1.914 \\
 20737 & 0.0442 & 1.74e-01 & 1.457 & 1.18e-03 & 1.995 & 1.26e-02 & 1.929 \\
 82433 & 0.0221 & 6.84e-02 & 1.350 & 2.95e-04 & 1.999 & 3.21e-03 & 1.972 \\
328705 & 0.0110 & 2.89e-02 & 1.245 & 7.38e-05 & 2.000 & 8.07e-04 & 1.990 \\
\hline
    \end{tabular}}
    \caption{Error history (error decay with respect to the number of degrees of freedom of the product discrete space, and experimental convergence rates) with respect to manufactured smooth solutions using the mixed-primal C0IP method with two different polynomial degrees.}
    \label{table:convergence}
\end{table}

\begin{figure}[t!]
    \centering
    \includegraphics[width=0.325\linewidth]{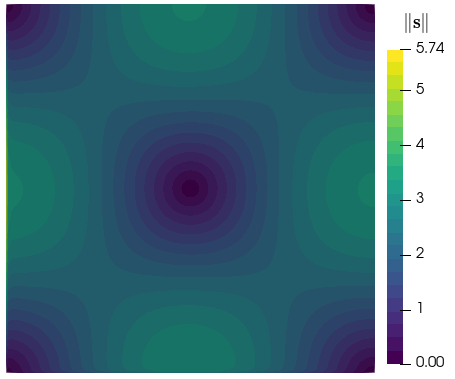}
    \includegraphics[width=0.325\linewidth]{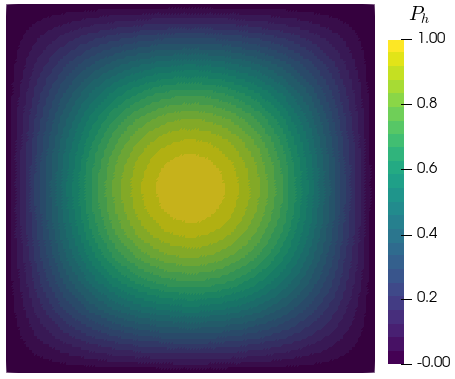}
    \includegraphics[width=0.325\linewidth]{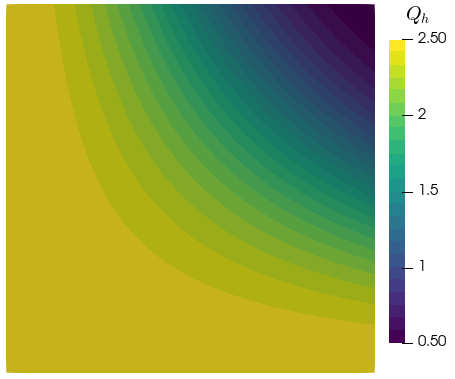}
    \caption{Approximate solutions for the convergence test, computed with the lowest-order finite element scheme and shown on a coarse mesh.}
    \label{fig:convergence}
\end{figure}

The results are presented in Table~\ref{table:convergence}, exhibiting optimal convergence for the density $O(h^{k+1})$, while the flux shows a sub-optimal error decay at $O(h^{k+1/2})$ (a sub-optimal convergence has been shown theoretically in \cite{neilan2019low} for axisymmetric Darcy equations approximated by lowest-order Raviart--Thomas elements). All unknowns exhibit an optimal convergence (in the sense that this is the approximability property of the chosen finite element spaces). Other tests (not shown here) confirm that if we do not use the C0IP strategy, while for the lowest-order case the density moment has optimal convergence, for the case $k=1$ we obtain only a $O(h)$ convergence. We also plot the components of the approximate solution in Figure~\ref{fig:convergence}. For this test, for every mesh refinement and polynomial degree the Newton--Raphson algorithm used for the nonlinear system has taken no more than four iterations to achieve the desired converge criterion. 

\section{Results}\label{sec:results}
In what follows we assume that the reaction boundary $\partial \Omega_R$ is defined as the boundary of the region
\begin{equation}\label{eq:exp_boundary}
     \Omega_{R} = \left\{(\pos{\eta}{1}, \pos{\eta}{2}) : r_1 \leq \sigmamax\text{exp}\left(\frac{-4\pi \Gamma r_2^3}{3}\right)\right\},
\end{equation}
where $\Gamma > 0$ is a parameter that controls the rate at which the height of the boundary decays as $r_2$ increases, and restrict our attention to the finite element domain
\begin{equation}\label{eq:finite_ele_domain}
    \Omega_{\text{FE}} = \left\{(r_1, r_2) : \sigmamax\text{exp}\left(\frac{-4\pi \Gamma r_2^3}{3}\right) \leq r_1 \leq 5,0 \leq r_2 \leq 5\right\}.
\end{equation}
This boundary is of particular interest because it was used by \textcite{our_first_paper} to implement a particle-based simulation of the system we are considering. Using Equation \cref{eq:closest_rate_old} we find that steady-state reaction rate for this domain can be approximated by
\begin{equation}\label{eq:exp_boundary_approx_flux}
K_1(c) = \frac{4\pi \hat{D}_1 \sigmamax c}{V_{\text{FE}}(c+\Gamma)},
\end{equation}
where $V_{\text{FE}} = (4/3)\pi (r_1^{\text{max}})^3$ and we have implicitly adopted the limit $r_2^{\text{max}} \rightarrow \infty$ which introduces a negligible error compared to integrating up until $r_2^{\text{max}} =5$. Although Kearney and Flegg briefly compared the reaction rate obtained from their simulation to $K_1(c)$, their results contained additional sources of error beyond those that arise from adopting this approximation making it impossible to isolate these contributions. To rectify this, in \cref{sec:finite_element_solutions} we examine the steady-state finite element solution for this particular choice of $\partial \Omega_R$ and compare the corresponding reaction rate with Equation \cref{eq:exp_boundary_approx_flux}. Then in \cref{sec:reactive_boundary_corrections} we demonstrate how to alter $\partial \Omega_R$ so that the observed reaction rate better replicates the desired reaction rate. In other words, we show how to correct for the errors incurred when adopting the approximation in Equation \cref{eq:exp_boundary_approx_flux}.\par

\subsection{Comparison of finite element and leading-order solutions}\label{sec:finite_element_solutions}

\begin{figure}[htp]
\begin{subfigure}{.5\textwidth}
  \centering
  \includegraphics[width=\textwidth]{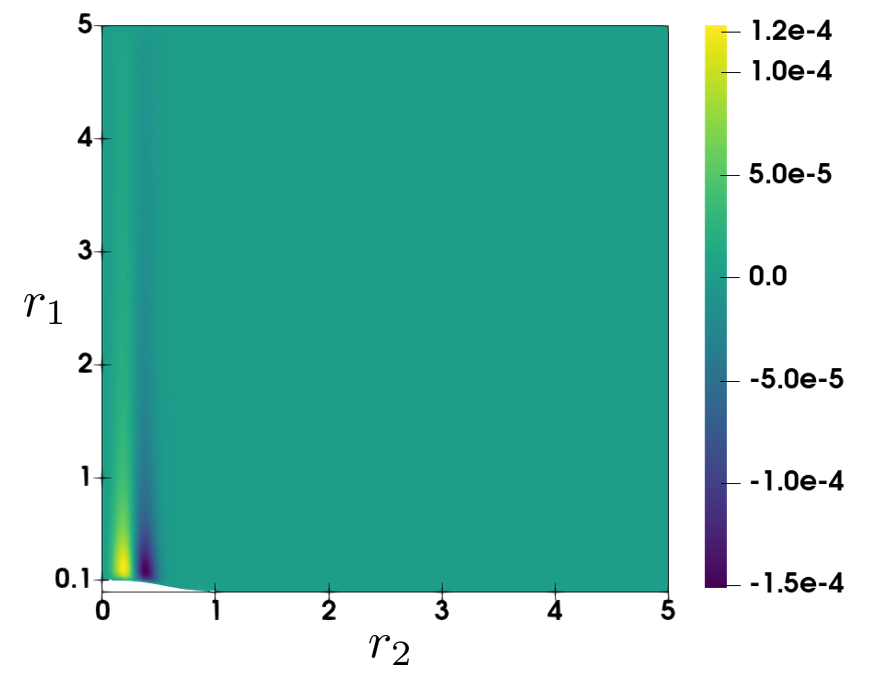}
  \caption{}
  \label{fig:dif_sigma0.1_gamma1_c10_nozoom}
\end{subfigure}
\begin{subfigure}{.5\textwidth}
  \centering
  \includegraphics[width=\textwidth]{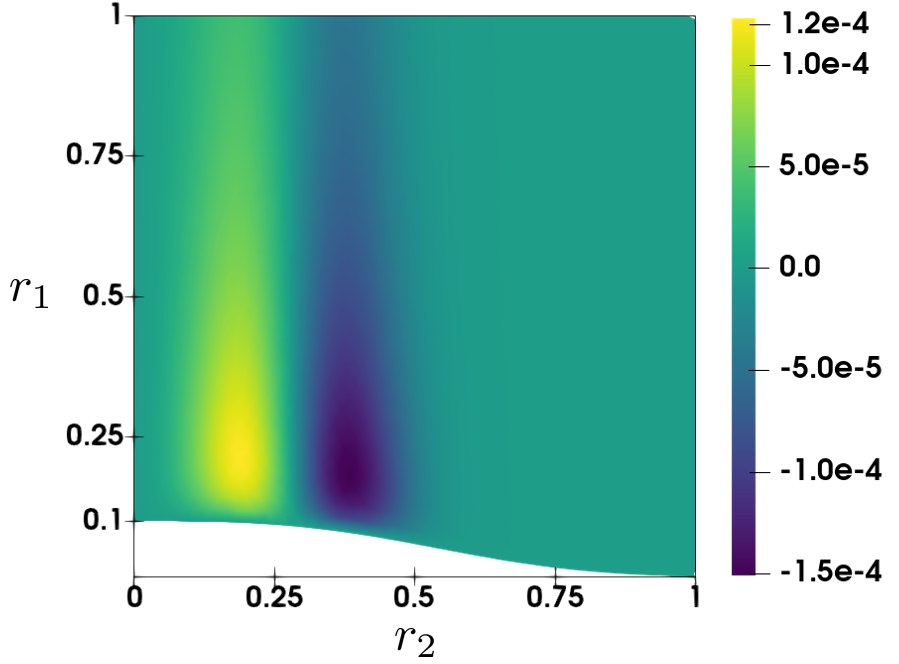}
  \caption{}
  \label{fig:dif_sigma0.1_gamma1_c10_zoom}
\end{subfigure}
\caption{The (scaled) difference $\Delta \steadyProb$, defined in Equation \cref{eq:scaled_difference}, between the steady-state finite element solution $\steadyProbfe$ and the leading-order solution $P_0(r_1,r_2)$ for $\hat{D}_1 = 2$, $\hat{D}_2 = 1.5$, $\sigmamax = 0.1$, $\Gamma = 1$, and $c=10$. The leading-order solution underestimates the probability, resulting in the yellow region where $\Delta \steadyProb > 0$, that the closest molecule of $C$ is near to the origin and instead favours larger values of $r_2$, which results in the dark blue region where $\Delta \steadyProb < 0$. This discrepancy arises because $P_0$ neglects the advection towards the origin in the $\pos{\hat{r}}{2}$ (horizontal) direction. The difference appears most pronounced near the bottom boundary $\partial \Omega_R$, defined in Equation \cref{eq:exp_boundary}, but this is because the solutions decay with increasing $r_1$. The surface plot in (a) shows the finite element domain from Equation \cref{eq:finite_ele_domain}, while the plot in (b) shows the subdomain $\left\{(r_1, r_2) : \sigmamax\text{exp}\left(\frac{-4\pi \Gamma r_2^3}{3}\right) \leq r_1 \leq 1,0 \leq r_2 \leq 1\right\}$.}
\label{fig:dif_sigma0.1_gamma1_c10}
\end{figure}

\begin{figure}[t!]
\begin{subfigure}{1\textwidth}
  \centering
  \includegraphics[width=0.75\textwidth]{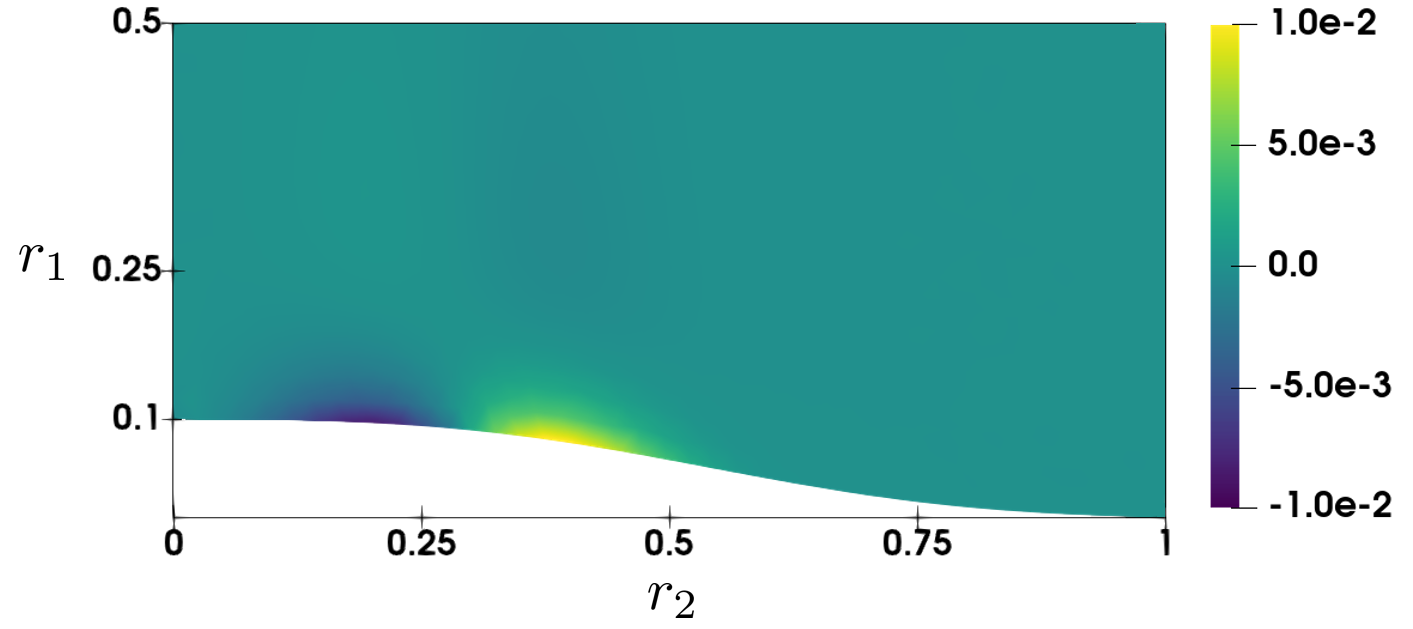}
  \caption{}
  \label{fig:dif_flux_r1_sigma0.1_gamma1_c10_zoom}
\end{subfigure}\\
\begin{subfigure}{1\textwidth}
  \centering
  \includegraphics[width=0.75\textwidth]{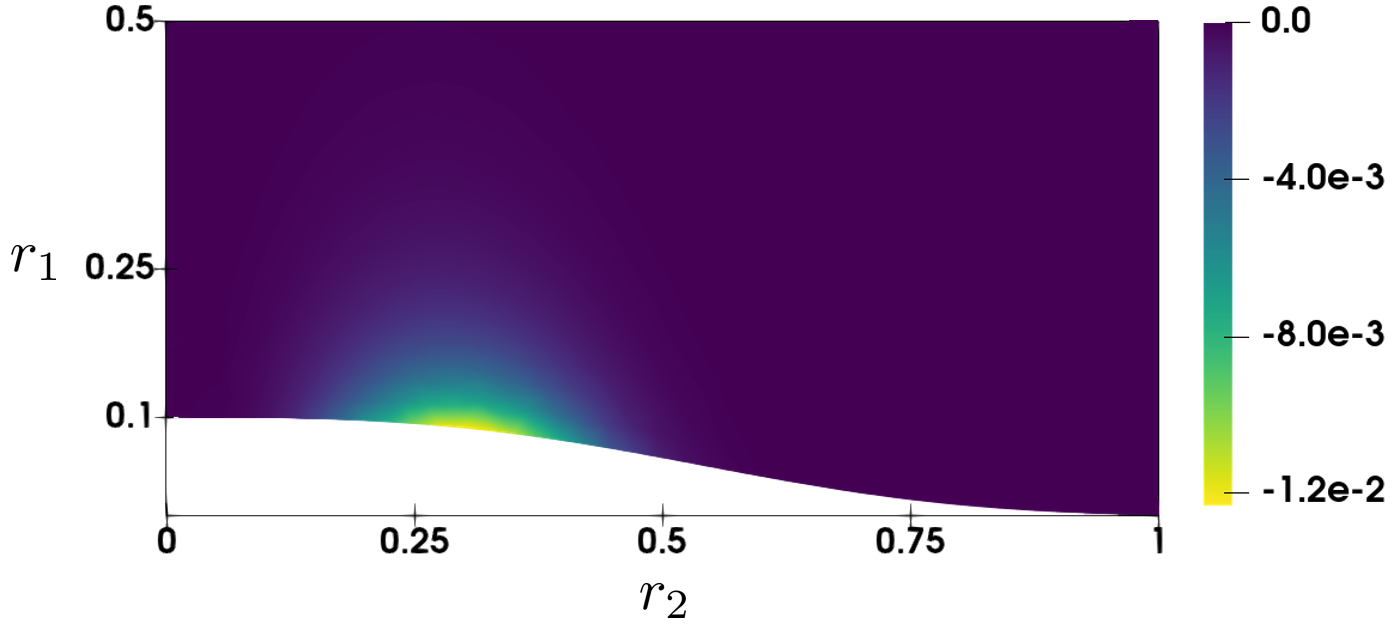}
  \caption{}
  \label{fig:dif_flux_r2_sigma0.1_gamma1_c10_zoom}
\end{subfigure}
\caption{The corrections to the flux of the 
 probability density that are neglected by the approximation $K_1(c)$ to the steady-state reaction rate defined in Equation \cref{eq:closest_rate_old}, for $\hat{D}_1 = 2$, $\hat{D}_2 = 1.5$, $\sigmamax = 0.1$, $\Gamma = 1$, and $c=10$. The surface plot in (a) shows the (scaled) correction $\Delta s_1(r_1, r_2)$ to the diffusive flux in the $\pos{\hat{r}}{1}$ (vertical) direction defined by Equation \cref{eq:r1_flux_correction}. The diffusive flux of the leading-order solution $P_0(r_1,r_2)$ underestimates the magnitude of the $\pos{\hat{r}}{1}$ component of the finite element solution $\sfe$ for values of $r_2$ less than approximately $0.25$, resulting in the dark blue region where $\Delta s_1(r_1, r_2) < 0$. In contrast, the magnitude of this flux is overestimated for larger values of $r_2$, resulting in the yellow region where $\Delta s_1(r_1, r_2) > 0$. Similarly, the surface plot in (b) shows the (scaled) correction $\Delta s_2(r_1, r_2)$ to the flux accounted for by $K_1(c)$ in the $\pos{\hat{r}}{2}$ (horizontal) direction defined by Equation \cref{eq:r2_flux_correction}. Since $K_1(c)$ entirely neglects the flux in this direction over $\partial \Omega_R$, as it is an $O(\sigmamax^2)$ contribution, this correction is simply the $\pos{\hat{r}}{2}$ component of the finite element solution $\sfe$. Both plots display the subdomain $\left\{(r_1, r_2) : 0 \leq r_1 \leq 1,0 \leq r_2 \leq 1\right\}$.}
\label{fig:dif_flux_sigma0.1_gamma1_c10}
\end{figure}

To isolate the $O(\sigmamax^2)$ corrections to the leading-order solution $P_0$, we consider the (scaled) difference $\Delta \steadyProb$ between $\steadyProbfe$ and Equation \cref{eq:P0_solution}
\begin{equation}\label{eq:scaled_difference}
    \Delta \steadyProb = 4\pi r_2^2\left(\steadyProbfe - \frac{c}{V_{\text{FE}}}\exp\left(\frac{-4\pi c r_2^3}{3}\right)\left(1-\frac{\sigmamax}{r_1}\text{exp}\left(\frac{-4\pi \Gamma r_2^3}{3}\right)\right)\right).
\end{equation}
The factor of $4\pi r_2^2$ is included since it is more convenient to work with the scaled probability density $4\pi r_2^2 \steadyProb$, which represents the probability to find the closest molecule of $C$ on the surface of a sphere of radius $r_2$ that is centred on the origin. \cref{fig:dif_sigma0.1_gamma1_c10} depicts $\Delta P$ in the case that $\hat{D}_1 = 2$, $\hat{D}_2 = 1.5$, $\sigmamax = 0.1$ and $\Gamma = 1$, while we set $c= 10$ since it results in a sufficiently concentrated density profile for the associated flux to be easily visible near $\partial \Omega_R$; see \cref{fig:dif_flux_sigma0.1_gamma1_c10}. The surface plot in \cref{fig:dif_sigma0.1_gamma1_c10_nozoom} displays the entirety of $\Omega_{\text{FE}}$ while \cref{fig:dif_sigma0.1_gamma1_c10_zoom} displays a subdomain that only extends to $r_1 = 1$ and $r_2 = 1$.\par

As expected, the solutions differ the most near $\partial \Omega_R$ while the difference decays as both $r_1$ and $r_2$ increase. Notice that while $\Delta P$ decays very rapidly to zero as $r_2$ increases, it decays much more slowly as $r_1$ increases and remains finite even very close to the boundary at $r_1 = r_1^{\text{max}}$; $\Delta P$ vanishes on this boundary since we enforce equality between the solutions there. This suggests that, similar to the $O(\sigmamax)$ contribution present in Equation \cref{eq:P0_solution}, the $O(\sigmamax^2)$ contributions to the solution decay exponentially as $r_2$ increases, but polynomially as $r_1$ increases. Moreover, we can see that $P_0$ underestimates the probability, resulting in the yellow region where $\Delta \steadyProb > 0$, that the closest molecule is close to the origin and instead favours larger values of $r_2$, resulting in the dark blue region where $\Delta \steadyProb < 0$. This is because the leading-order solution results from considering the limit $\sigmamax^2 \rightarrow 0$ in Equation \cref{eq:perturbation_pde}, which amounts to ignoring both the diffusive motion in the $\pos{\hat{r}}{2}$ direction and, more importantly, the advective motion towards the origin in this direction. \par

Similarly, when approximating the steady-state reaction rate using $K_1(c)$ from Equation \cref{eq:closest_rate_old}, we only consider the $O(\sigmamax)$ contributions to the flux of the probability density over $\partial \Omega_R$. This effectively ignores the contribution to this reaction rate from the flux in the $\pos{\hat{r}}{2}$ direction as well as an $O(\sigmamax^2)$ correction to the contribution in the $\pos{\hat{r}}{1}$ direction, which arises from the correction of the same order to $P_0$. The (scaled) $O(\sigmamax^2)$ correction to the diffusive flux in the $\pos{\hat{r}}{1}$ direction is given by
\begin{equation}\label{eq:r1_flux_correction}
    \Delta s_1(r_1, r_2) = 4\pi r_2^2\hat{D}_1 \left(\sfe \cdot \pos{\hat{r}}{1} + \nabla P_0(r_1, r_2) \cdot \pos{\hat{r}}{1} \right), 
\end{equation}
and is shown in \cref{fig:dif_flux_r1_sigma0.1_gamma1_c10_zoom}, where in this case
\begin{equation}
    \nabla P_0(r_1, r_2)\cdot \pos{\hat{r}}{1} = \frac{ c \sigmamax}{V_{\text{FE}}r_1^2}\exp\left(\frac{-4\pi (c+\Gamma) r_2^3}{3}\right).
\end{equation} Since $K_1(c)$ ignores the flux in the $\pos{\hat{r}}{2}$ direction entirely, the (scaled) $O(\sigmamax^2)$ correction to this component is simply the total flux in this direction
\begin{equation}\label{eq:r2_flux_correction}
   \Delta s_2(r_1, r_2) = 4\pi r_2^2 \hat{D}_2\sfe \cdot \pos{\hat{r}}{2},
\end{equation}
and is shown in \cref{fig:dif_flux_r2_sigma0.1_gamma1_c10_zoom}. In each case, we have restricted our attention to a subdomain of $\Omega_{\text{FE}}$ that only extends to $r_1 = 0.5$ and $r_2 = 1$, since the flux is negligible outside this subdomain. The behaviour of $\Delta s_1$ mimics that of $\Delta P$ in that the flux of the leading-order solution underestimates the magnitude of the diffusive flux in the $\pos{\hat{r}}{1}$ direction for small values of $r_2$, resulting in the dark blue region where $\Delta s_1(r_1, r_2) < 0$, while overestimating it at larger values, resulting in the yellow region where $\Delta s_1(r_1, r_2) > 0$. In addition, we can see that there will also be a correction in the same vicinity due to $\Delta s_2$. Notice that both components of the flux point in the opposite direction to the respective radial unit vectors $\pos{\hat{r}}{1}$ and $\pos{\hat{r}}{2}$. \par

To better contextualise these corrections we now consider the total flux of $P_h$ over $\partial \Omega_R$
\begin{equation}\label{eq:total_FE_flux_s}
   K_\text{FE}(c) = \int_0^{r_2^{\text{max}}} D\sfe \cdot \pos{\hat{n}}{\text{out}}(4\pi f(r_2)r_2)^2 \,\mathrm{d}r_2, 
\end{equation}
which serves as a proxy for the exact reaction rate $K(c)$, although we will continue to use the subscript 'FE' to distinguish this approximate rate. We are concerned with how the corrections to the flux shown in \cref{fig:dif_flux_sigma0.1_gamma1_c10} affect this rate, and to quantify these contributions, we define three relative errors,

\begin{subequations}
\begin{equation}\label{eq:total_rel_error}
    \Delta K(c) = \frac{K_{\text{FE}}(c)-K_1(c)}{K_1(c)},
\end{equation} 

\begin{equation}\label{eq:r1_rel_error}
    \Delta K^1(c) = \frac{1}{K_1(c)}\left(\int_0^{r_2^{\text{max}}} D\sfe \cdot \pos{\hat{n}}{1}(4\pi f(r_2)r_2)^2 \,\mathrm{d}r_2-K_1(c)\right),
\end{equation}
and
\begin{equation}\label{eq:r2_rel_error}
    \Delta K^2(c) = \frac{1}{K_1(c)}\int_0^{r_2^{\text{max}}} D\sfe \cdot \pos{\hat{n}}{2}(4\pi f(r_2)r_2)^2 \,\mathrm{d}r_2,
\end{equation}\end{subequations}
where $\pos{\hat{n}}{i} = (\pos{\hat{n}}{\text{out}}\cdot \pos{\hat{r}}{i})\pos{\hat{r}}{i}$. Since $K_1(c)$ only considers the diffusive flux of $P_0$ over $\partial \Omega_R$ in the $\pos{\hat{r}}{1}$ direction, $\Delta K^1(c)$ essentially measures the relative error in the reaction rate that arises from neglecting $O(\sigmamax^2)$ corrections to the flux of $P$ in this direction. Similarly, $\Delta K^2(c)$ measures the relative error that arises from neglecting the diffusive and advective flux of $P$ over $\partial \Omega_R$ in the $\pos{\hat{r}}{2}$ direction.\par 

To explore how these errors behave as the shape of $\partial \Omega_R$ changes, we fix $c = 1$, $\hat{D}_1 = 2$, $\hat{D}_2 = 1.5$ and conduct two tests. In the first, we hold $\sigmamax = 0.1$ constant and increase $\Gamma$ from $0.25$ to $3$ in intervals of $0.25$, while in the second we instead hold $\Gamma = 1$ constant and increase $\sigmamax$ from $0.025$ to $0.25$ in intervals of $0.025$.\par 

\begin{figure}[t!]
\centering
  \includegraphics[width=1\textwidth]{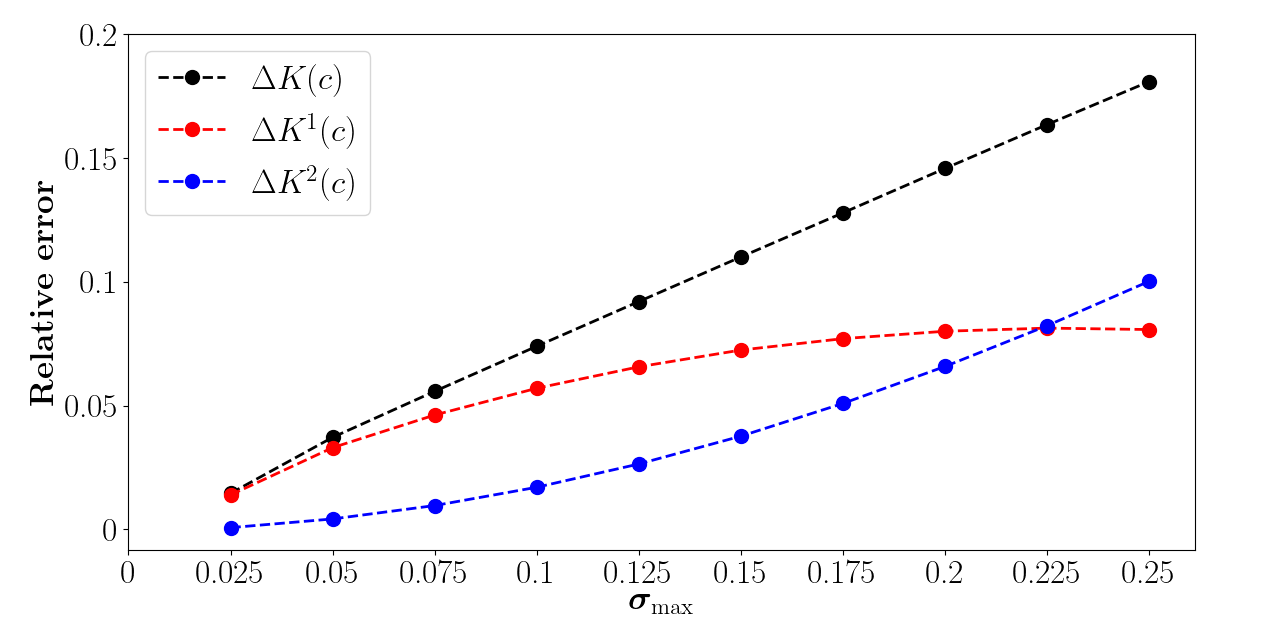}
  \label{fig:total_flux_sigma_test}

\caption{The relative error between the steady-state reaction rate $K_{\text{FE}}(c)$ obtained from the finite element solution via Equation \cref{eq:total_FE_flux_s}, and the approximate reaction rate $K_1(c)$, given by Equation \cref{eq:closest_rate_old}, calculated for $c = 1$, $\hat{D}_1 = 2$, $\hat{D}_2 = 1.5$, $\Gamma = 1$, while $\sigmamax$ is increased from $0.025$ to $0.25$ in intervals of $0.025$. The steady-state reaction rate is defined as the total flux of the probability density over the boundary $\partial \Omega_R$, defined by Equation \cref{eq:exp_boundary}. The black line shows the total relative error $\Delta K(c)$, defined in Equation \cref{eq:total_rel_error}, between $K_{\text{FE}}(c)$ and $K_1(c)$. This error arises because $K_1(c)$ neglects $O(\sigmamax^2)$ contributions to the flux over $\partial \Omega_R$. Since there are contributions of this order over $\partial \Omega_R$ in both the $\pos{\hat{r}}{1}$ and $\pos{\hat{r}}{2}$ directions, we quantify them separately using the relative errors $\Delta K^1(c)$ and $\Delta K^2(c)$, defined in equations \cref{eq:r1_rel_error} and \cref{eq:r2_rel_error}, and shown by the red, and blue lines respectively. That is, $\Delta K^1(c)$ and $\Delta K^2(c)$ represent the relative error that arises from ignoring $O(\sigmamax^2)$ contributions to the flux over $\partial \Omega_R$ in the $\pos{\hat{r}}{1}$ and $\pos{\hat{r}}{2}$ directions respectively. We can see that while $\Delta K(c)$ and $\Delta K^2(c)$ increase monotonically as $\sigmamax$ increases, $\Delta K^1(c)$ only increases until $\sigmamax = 0.225$ after which it begins to decrease. This behaviour arises because as $\sigmamax$ increases the magnitude of the $\pos{\hat{r}}{1}$ ($\pos{\hat{r}}{2}$) component of the outward normal vector to $\partial \Omega_R$ decreases (increases), thus decreasing (increasing) the flux over $\partial \Omega_R$ in the $\pos{\hat{r}}{1}$ ($\pos{\hat{r}}{2}$) direction. The finite element solution used to calculate $K_{\text{FE}}(c)$ was evaluated over the domain $\Omega_{\text{FE}} = \left\{(r_1, r_2) : \sigmamax\text{exp}\left(\frac{-4\pi \Gamma r_2^3}{3}\right) \leq r_1 \leq 5,0 \leq r_2 \leq 5\right\}$.}
\label{fig:flux_sigma_test}
\end{figure}

The black, red and blue lines in \cref{fig:flux_sigma_test} depict the relative errors $\Delta K(c)$, $\Delta K^1(c)$ and $\Delta K^2(c)$ respectively, in the case where $\sigmamax$ is varied and $\Gamma$ is held constant. In this case, $\Delta K(c)$ and $\Delta K^2(c)$ increase monotonically as $\sigmamax$ increases, but $\Delta K^1(c)$ only increases until $\sigmamax = 0.225$ after which it begins to decrease. Moreover, at the same point $\Delta K^2(c)$ becomes greater than $\Delta K^1(c)$ for the first time. This behaviour arises because the magnitude of the $\pos{\hat{r}}{2}$ component of the normal to the boundary, $\pos{\hat{n}}{\text{out}}$, increases with increasing $\sigmamax$, while the magnitude of the $\pos{\hat{r}}{1}$ component decreases (see Equation \cref{eq:normal_vector}). Therefore, the flux over $\partial \Omega_R$ in the $\pos{\hat{r}}{2}$ direction increases, while the flux over this boundary in the $\pos{\hat{r}}{1}$ direction decreases. Finally, recalling that $K_1(c)$ is $O(\sigmamax)$, we note that $\Delta K(c)$ growing linearly with $\sigmamax$ indicates that the corrections to $K_1(c)$ are $O(\sigmamax^2)$ as expected.\par

\begin{figure}[t!]
\centering
\includegraphics[width=1\textwidth]{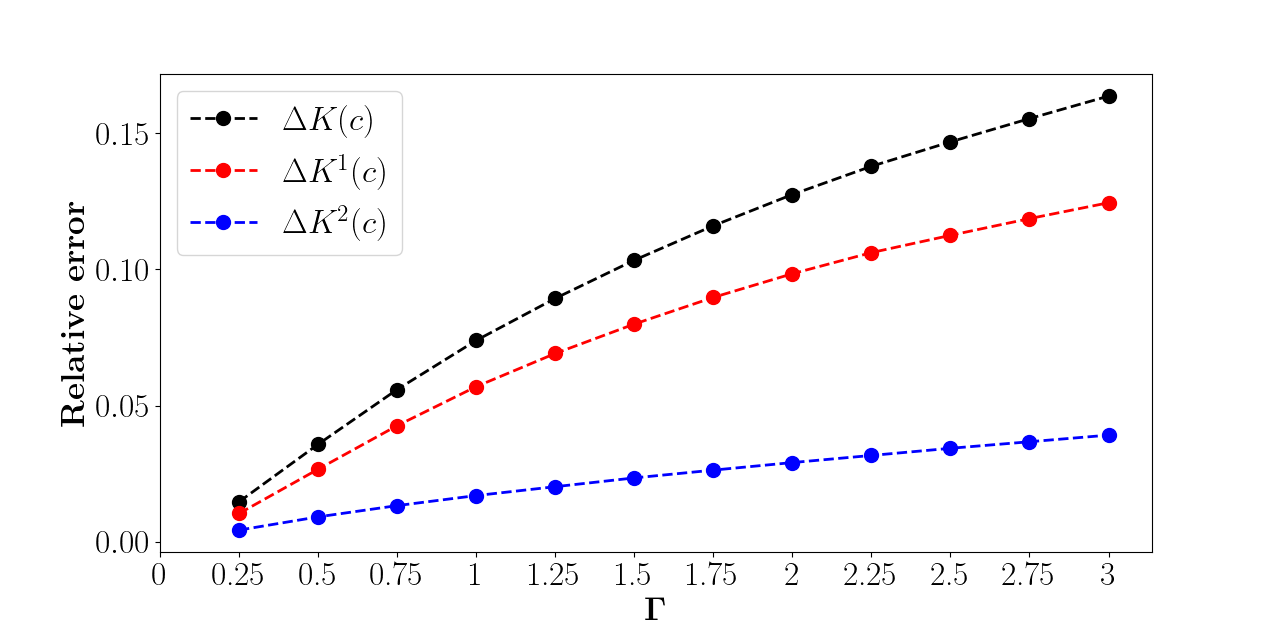}
\caption{The relative error between the steady-state reaction rate $K_{\text{FE}}(c)$ obtained from the finite element solution via Equation \cref{eq:total_FE_flux_s}, and the approximate reaction rate $K_1(c)$, given by Equation \cref{eq:closest_rate_old}, calculated for $c = 1$, $\hat{D}_1 = 2$, $\hat{D}_2 = 1.5$, and $\sigmamax = 0.1$, while $\Gamma$ is increased from $0.25$ to $3$ in intervals of $0.25$. The steady-state reaction rate is defined as the total flux of the probability density over the boundary $\partial \Omega_R$, defined by Equation \cref{eq:exp_boundary}. The black line shows the total relative error $\Delta K(c)$, defined in Equation \cref{eq:total_rel_error}, between $K_{\text{FE}}(c)$ and $K_1(c)$. This error arises because $K_1(c)$ neglects $O(\sigmamax^2)$ contributions to the flux over $\partial \Omega_R$. Since there are contributions of this order over $\partial \Omega_R$ in both the $\pos{\hat{r}}{1}$ and $\pos{\hat{r}}{2}$ directions, we quantify them separately using the relative errors $\Delta K^1(c)$ and $\Delta K^2(c)$, defined in equations \cref{eq:r1_rel_error} and \cref{eq:r2_rel_error}, and shown by the red, and blue lines respectively. That is, $\Delta K^1(c)$ and $\Delta K^2(c)$ represent the relative error that arises from ignoring $O(\sigmamax^2)$ contributions to the flux over $\partial \Omega_R$ in the $\pos{\hat{r}}{1}$ and $\pos{\hat{r}}{2}$ directions respectively. The relative errors increase monotonically, but the rate of this increase decreases with increasing $\Gamma$. The approximate rate $K_1(c)$ is exact when $\partial \Omega_R$ is constant, that is, when the outward unit normal vector to this boundary is $\pos{\hat{n}}{\text{out}} = -\pos{\hat{r}}{1}$, and so the relative errors will be large when $\pos{\hat{n}}{\text{out}}$ is poorly approximated by $-\pos{\hat{r}}{1}$, and will be small when this is a good approximation. With this in mind, the observed behaviour arises because for any particular $r_2$ value, which defines a particular point on $\partial \Omega_R$, the magnitude of the $\pos{\hat{r}}{2}$ ($\pos{\hat{r}}{1}$) component of $\pos{\hat{n}}{\text{out}}$ will initially increase (decrease), eventually reach a local maximum (minimum) and then begin to decrease (increase) as $\Gamma$ is increased. The finite element solution used to calculate $K_{\text{FE}}(c)$ was evaluated over the domain $\Omega_{\text{FE}} = \left\{(r_1, r_2) : \sigmamax\text{exp}\left(\frac{-4\pi \Gamma r_2^3}{3}\right) \leq r_1 \leq 5,0 \leq r_2 \leq 5\right\}$.}
\label{fig:flux_gamma_test}
\end{figure}

Similarly, \cref{fig:flux_gamma_test} displays the relative errors $\Delta K(c)$, $\Delta K^1(c)$ and $\Delta K^2(c)$ as a black, red, and blue line, respectively, in the case where $\Gamma$ is varied and $\sigmamax$ is constant. The three relative errors increase monotonically, but the rate of increase appears to decrease as $\Gamma$ increases. The height of $\partial \Omega_R$ decays more rapidly as $\Gamma$ increases, but the relationship between this parameter and $\pos{\hat{n}}{\text{out}}$ is less straightforward. For any particular value of $r_2$, which defines a particular point on $\partial \Omega_R$, the magnitude of the $\pos{\hat{r}}{2}$ ($\pos{\hat{r}}{1}$) component of $\pos{\hat{n}}{\text{out}}$ will initially increase (decrease) with increasing $\Gamma$ until it reaches a local maximum (minimum) after which it begins to decrease (increase) as $\Gamma$ increases further. Moreover, the value of $\Gamma$ at which this stationary point occurs decreases as $r_2$ increases. Therefore, for a fixed concentration and hence density profile, the approximation $\pos{\hat{n}}{\text{out}} = -\pos{\hat{r}}{1}$ initially becomes worse as $\Gamma$ increases from zero and the relative errors all increase accordingly. However, if $\Gamma$ continues to increase, eventually the magnitude of the $\pos{\hat{r}}{2}$ component of $\pos{\hat{n}}{\text{out}}$ in the region where the density is the most concentrated will begin to decay and $\pos{\hat{n}}{\text{out}}$ begins to converge to $-\pos{\hat{r}}{1}$ and so $K_\text{FE}(c)$ begins to converge to $K_1(c)$. Therefore, the relative errors will initially increase, plateau and eventually begin to decrease again once $\Gamma$ becomes large enough.
\par

\subsection{Reaction boundary corrections}\label{sec:reactive_boundary_corrections}
To reproduce realistic chemical kinetics, particle-based simulations must accurately reproduce known reaction rates. In other words, modellers typically begin with a desired reaction rate, which we denote $K_D(c)$, and then determine what reaction conditions should be imposed on the system to reproduce this rate. The expression for $K_1(c)$ given in Equation \cref{eq:closest_rate_old} is advantageous since it can be easily solved for the function $f(r_2)$, which defines $\partial \Omega_R$, via an inverse Laplace transform \cite{our_first_paper}. That is, it allows us to derive an analytic expression for the reaction boundary $\partial \Omega_R$ that (approximately) matches a given rate $K_1(c)$. Thus, the obvious approach is to replace $K_1(c)$ with the desired rate $K_D(c)$ and determine the corresponding reaction boundary. However, our results demonstrate that deriving $\partial \Omega_R$ in this way will result in a boundary that does not reproduce $K_D(c)$ exactly. Instead, the simulated reaction rate (assuming that all other sources of error are negligible) will be given by Equation \cref{eq:total_flux_s}. This is problematic since small changes in the reaction rate can drastically change the behaviour of biochemical systems (for example, in systems that exhibit bistability and hysteresis \cite{bistability_hysteresis, wilhelm2009smallest}) and, for the specific form of $\partial \Omega_R$ considered here, the relative error between $K_1(c)$---or in this case $K_D(c)$---and $K_\text{FE}(c)$ becomes large as $\sigmamax$ or $\Gamma$ increases. Therefore, we would like to correct $\partial \Omega_R$ so that the reaction rate that $\partial \Omega_R$ actually reproduces, which we quantify using $K_\text{FE}(c)$, more closely resembles the desired reaction rate $K_D(c)$. Moreover, since during a simulation $\partial \Omega_R$ must be evaluated each time there is potential for a reaction to occur, we would like to retain an analytic expression for this boundary that is efficient to evaluate and seek to alter the parameters that define $\partial \Omega_R$ rather than its functional form. This amounts to assuming that the functional form of $K_\text{FE}(c)$ is well approximated by that of $K_1(c)$.\par 

\begin{figure}[t!]
\begin{subfigure}{1\textwidth}
  \centering
  \includegraphics[width=0.75\textwidth]{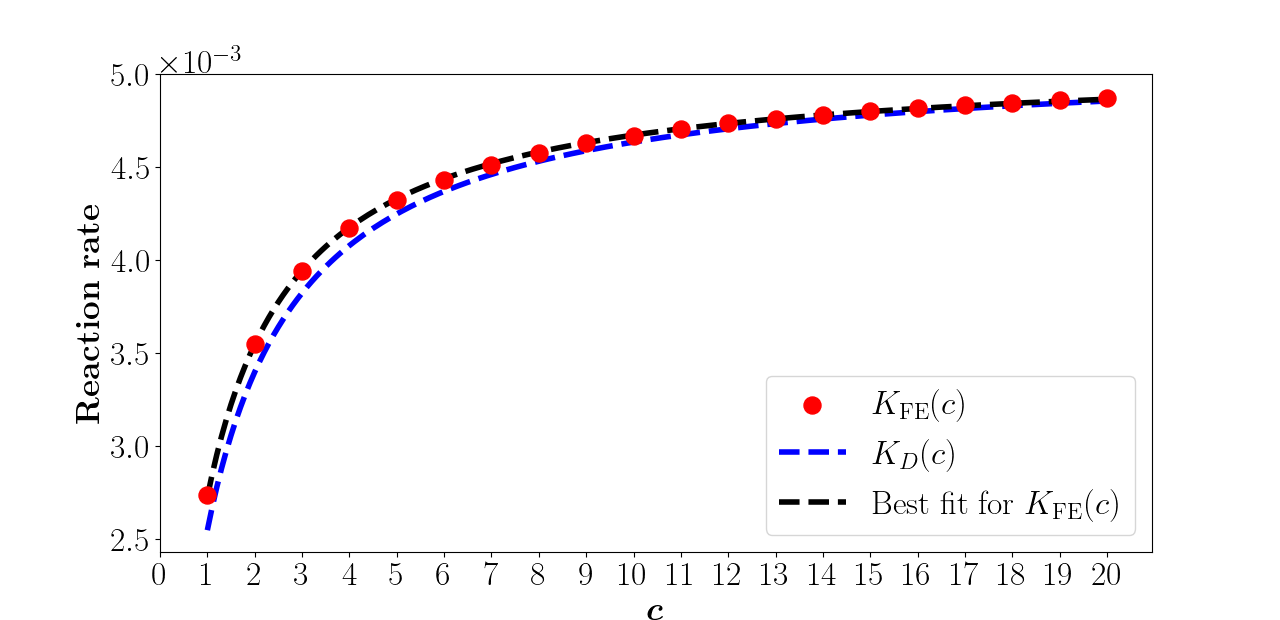}
  \caption{}
  \label{fig:uncorrected_rates}
\end{subfigure}\\
\begin{subfigure}{1\textwidth}
  \centering
  \includegraphics[width=0.75\textwidth]{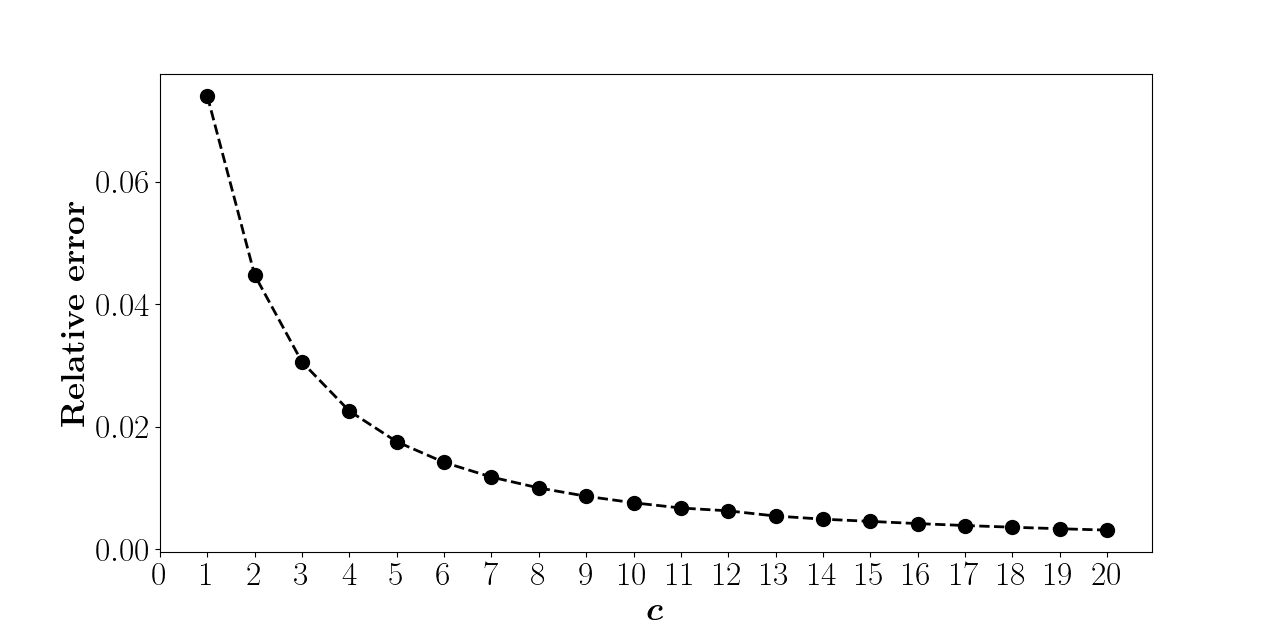}
  \caption{}
  \label{fig:uncorrected_relative_error}
\end{subfigure}
\caption{The steady-state reaction rate $K_\text{FE}(c)$ obtained from the finite element solution via Equation \cref{eq:total_FE_flux_s}, for $\sigmamax = 0.1$, $\Gamma = 1$, $\hat{D}_1 = 2$, and $\hat{D}_2 = 1.5$, while the concentration $c$ is increased from $1$ to $20$ in increments of 1. The plot in (a) compares $K_\text{FE}(c)$, shown by the red points, to the desired reaction rate $K_D(c)$ shown by the blue line and defined in Equation \cref{eq:desired_reaction_rate} with $\sigmainit = 0.1$ and $\gammainit = 1$. The black line shows the line of best fit for $K_\text{FE}(c)$ when assuming the model in Equation \cref{eq:leading_order_model} and is defined by the parameters $\sigmaeff = 0.1$ and $\gammaeff = 0.86$. The plot in (b) displays the relative error $\Delta K_D(c)$, defined in Equation \cref{eq:desired_rate_rel_error}, between $K_\text{FE}(c)$ and $K_D(c)$ for the same set of parameters and shows that $\Delta K_D(c)$ decreases monotonically as $c$ is increased.} 
\label{fig:uncorrected_concentration_test}
\end{figure}

To demonstrate how to correct $\partial \Omega_R$ in this manner we suppose that the desired reaction rate is given by \begin{equation}\label{eq:desired_reaction_rate}
    K_D(c) = \frac{\sigmainit c}{V_{\text{FE}}(c + \gammainit)},
\end{equation}
with $\sigmainit = 0.1$ and $\gammainit = 1$. Using Equation \cref{eq:closest_rate_old} we find that the boundary defined in Equation \cref{eq:exp_boundary} (approximately) matches this reaction rate. Therefore, we adopt this functional form for $\partial \Omega_R$ and set $\sigmamax = 0.1$ and $\Gamma = 1$. In \cref{fig:uncorrected_concentration_test} we compare $K_D(c)$ and $K_{\text{FE}}(c)$ for this initial boundary in the case that $\hat{D}_1 = 2$, $\hat{D}_2 = 1.5$ and $c$ is increased from $1$ to $20$ in increments of $1$. \cref{fig:uncorrected_rates} depicts $K_D(c)$ and $K_{\text{FE}}(c)$ using a blue line and the red points, respectively. The corresponding relative error,
\begin{equation}\label{eq:desired_rate_rel_error}
    \Delta K_D(c) = \frac{K_{\text{FE}} - K_D(c)}{K_D(c)},
\end{equation}
is shown in \cref{fig:uncorrected_relative_error}
and decreases monotonically as $c$ increases. This occurs because as $c$ increases, the probability of finding the closest molecule of $C$ close to the origin (at small $r_2$) increases and the probability density profile narrows, which reduces the influence the curvature of $\partial \Omega_R$ has on the reaction rate. In other words, as $P_h$ becomes more localised the relevant part of $\partial \Omega_R$ (in the region where $4 \pi r_2^2\steadyProb > 0$) is better approximated by a constant. This means that $K_{\text{FE}}(c)$ is well approximated by $K_1(c)$ and therefore matches $K_D(c)$ more closely since in this case $K_D(c) = K_1(c)$.\par

To better quantify the difference between $K_{\text{FE}}(c)$ and $K_D(c)$ across the range of concentrations considered we assume $K_{\text{FE}}(c)$ is well approximated by the functional form of $K_D(c)$, and fit the model
\begin{equation}\label{eq:leading_order_model}
    K_{\text{FE}}(c) = \frac{\sigmaeff c}{V_{\text{FE}}(c + \gammaeff)},
\end{equation}
to determine the `effective' parameters $\sigmaeff$ and $\gammaeff$. If $\partial \Omega_R$ reproduced the desired reaction rate exactly, then $\sigmaeff$ and $\gammaeff$ would be identical to $\sigmainit$ and $\gammainit$ respectively. However, we already know that this is not the case here and the resulting line of best fit is defined by $\sigmaeff = 0.1$ and $\gammaeff = 0.86$ and is shown by the black dotted line in \cref{fig:uncorrected_rates}. Through this fit we can see that $\sigmainit$ is well approximated by $\sigmaeff$, but there is a large discrepancy between $\gammainit$ and $\gammaeff$ with a corresponding relative error of $14\%$.\par

Our goal now is to alter $\partial \Omega_R$, by varying $\sigmamax$ and $\Gamma$, so that the error between $K_{\text{FE}}(c)$ and $K_D(c)$ is minimised. It is important to emphasise that $K_D(c)$ is determined by the modeller, and that our goal is to determine a reaction boundary that reproduces this rate. In contrast, $K_1(c)$ and $K_{\text{FE}}(c)$ both approximate the reaction rate associated with the current reaction boundary $\partial \Omega_R$. The former does so through Equation \cref{eq:closest_rate_old}, while the latter is connected to $\partial \Omega_R$ through Equation \cref{eq:total_FE_flux_s} and serves as a proxy for the exact reaction rate $K(c)$. Crucially, altering $\partial \Omega_R$ changes $K_1(c)$ and $K_{\text{FE}}(c)$, but not $K_D(c)$.\par

Since the error between $K_D(c)$ and $K_{\text{FE}}(c)$ is dependent on the concentration, we use the Levenberg-Marquardt algorithm \cite{Levenberg, Marquardt, non-linear_least_squares_book} to minimise the sum of the squares of the residuals between these quantities. We define the residual for a given concentration as
\begin{equation}
    r(c) = K_D(c) - K_{\text{FE}}(c),
\end{equation}
and recall that initially $\sigmamax = 0.1$ and $\Gamma = 1$ define $\partial \Omega_R$ through Equation \cref{eq:exp_boundary}. In addition, we use Equation \cref{eq:exp_boundary_approx_flux} to derive an analytic approximation of the Jacobian for $K_{\text{FE}}(c)$. Each iteration of the algorithm we: 
\begin{enumerate}
    \item construct a finite element mesh $\mathcal{T}_h$ according to the current values of $\sigmamax$ and $\Gamma$,
    \item compute $K_{\text{FE}}(c)$, increasing $c$ from $1$ to $20$ in intervals of size $1$,
    \item update the values of the parameters $\sigmamax$ and $\Gamma$, thus defining the new boundary $\partial \Omega_R$ and Jacobian for the next iteration.
\end{enumerate}
This process is continued until the sum of the squares of the residuals or the values of $\sigmamax$ and $\Gamma$ do not change significantly between successive iterations. \par 

\begin{figure}[t!]
\begin{subfigure}{1\textwidth}
  \centering
  \includegraphics[width=0.75\textwidth]{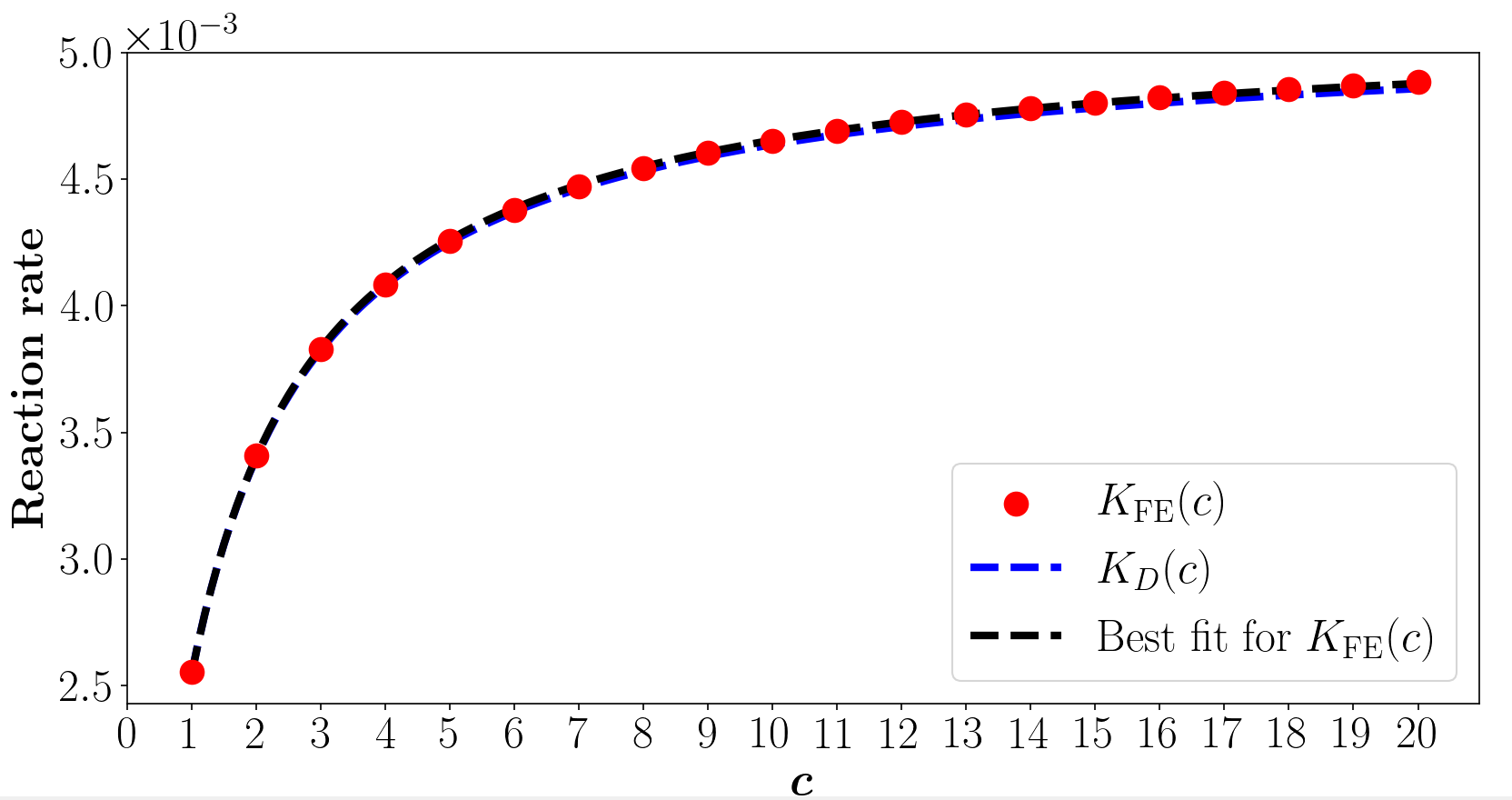}
  \caption{}
  \label{fig:corrected_rates}
\end{subfigure}\\
\begin{subfigure}{1\textwidth}
  \centering
  \includegraphics[width=0.75\textwidth]{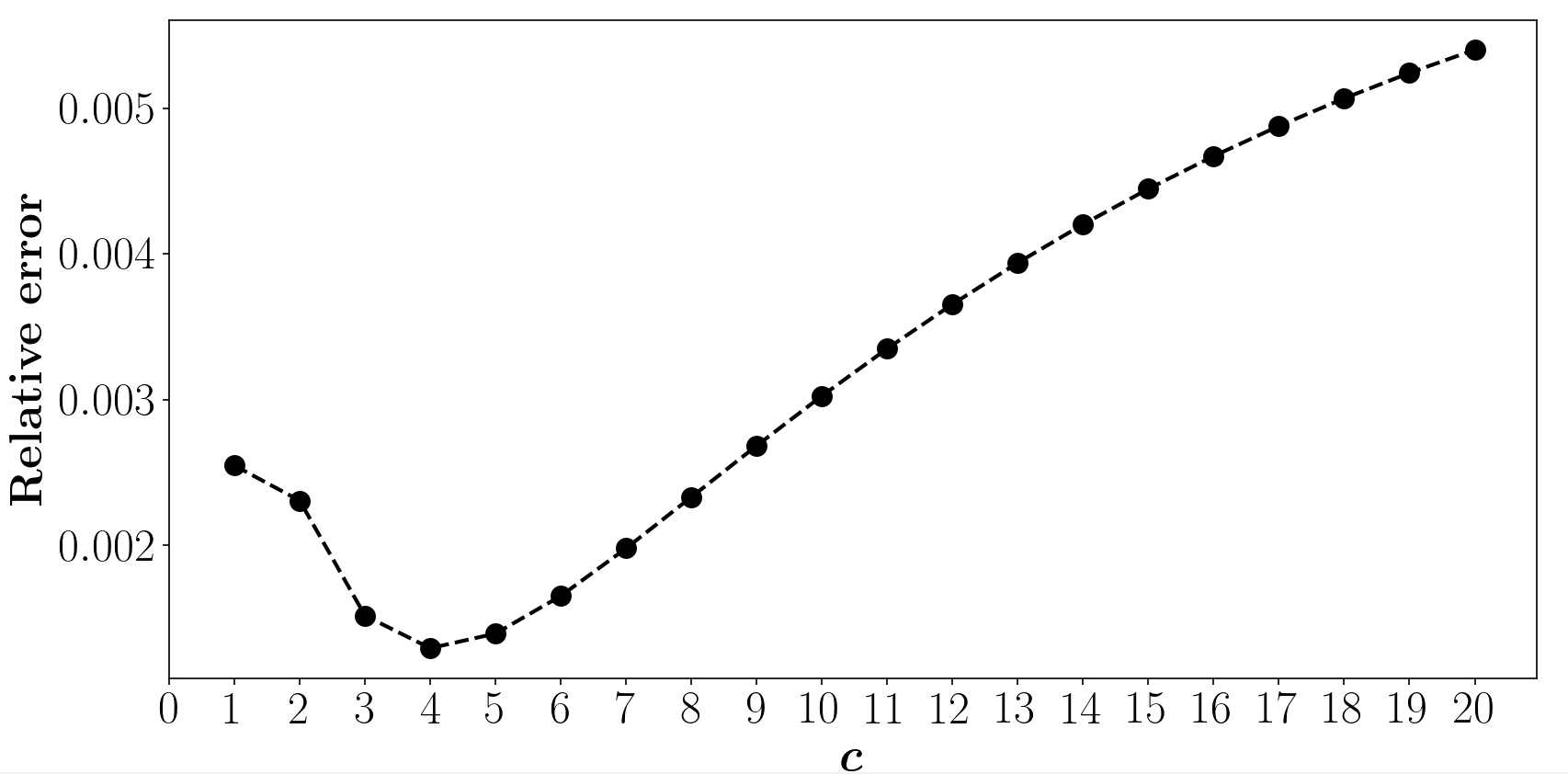}
  \caption{}
  \label{fig:corrected_relative_error}
\end{subfigure}
\caption{The steady-state reaction rate $K_\text{FE}(c)$ obtained from the finite element solution via Equation \cref{eq:total_FE_flux_s}, for the corrected boundary defined by Equation \cref{eq:exp_boundary} with $\sigmamax = 0.1$ and $\Gamma = 1.19$. To enable comparison with \cref{fig:uncorrected_concentration_test} we again consider $\hat{D}_1 = 2$ and $\hat{D}_2 = 1.5$, while the concentration $c$ is increased from $1$ to $20$ in increments of 1. The plot in (a) compares $K_\text{FE}(c)$, shown by the red points, to the desired reaction rate $K_D(c)$ shown by the blue line and defined in Equation \cref{eq:desired_reaction_rate} with $\sigmainit = 0.1$ and $\gammainit = 1$. The black line shows the line of best fit for $K_\text{FE}(c)$ when assuming the model in Equation \cref{eq:leading_order_model} and is defined by the parameters $\sigmaeff = 0.1$ and $\gammaeff = 1.01$. The plot in (b) displays the relative error $\Delta K_D(c)$, defined in Equation \cref{eq:desired_rate_rel_error}, between $K_\text{FE}(c)$ and $K_D(c)$ for the same set of parameters. When compared to \cref{fig:uncorrected_relative_error} we can see that the corrected boundary significantly reduces the relative error for $c \leq 15$. However, now the relative error tends to increase with increasing $c$ and is worse than it was previously for $c > 15$. This occurs because for the uncorrected boundary, defined instead by $\sigmamax = 0.1$ and $\Gamma = 1$, 
$K_\text{FE}(c)$ converges to $K_D(c)$ in the limit that $c \rightarrow \infty$ and so altering the boundary will increase rather than decrease the relative error if $c$ is sufficiently large.}
\label{fig:corrected_concentration_test}
\end{figure}

In our case, the optimisation yields $\sigmamax = 0.10$ and $\Gamma = 1.19$ and we compare $K_{\text{FE}}(c)$ for this boundary to $K_D(c)$ in \cref{fig:corrected_concentration_test}. The blue line in \cref{fig:corrected_rates} represents $K_D(c)$, the red points denote $K_{\text{FE}}(c)$ calculated for the boundary in Equation \cref{eq:exp_boundary} defined by the optimised values of $\sigmamax$ and $\Gamma$, and the black line is the line of best fit through these points for the model in Equation \cref{eq:leading_order_model}. The new line of best fit is defined by $\sigmaeff = 0.1$ and $\gammaeff = 1.01$, which means that again $\sigmainit$ is well approximated by $\sigmaeff$, while the relative error between $\gammainit$ and $\gammaeff$ has reduced in magnitude from $14\%$ to just $1\%$. Similarly, \cref{fig:corrected_relative_error} shows a reduction in the relative error between $K_{\text{FE}}(c)$ and $K_D(c)$ for $c \leq 15$ while there is a slight increase in the error for $c > 15$. Finally, the sum of the squares of the residuals has been significantly reduced from $2.5 \times 10^{-2}$ to $1 \times 10^{-3}$.

\section{Discussion}
The leading-order solution $P_0$ to the steady-state probability density in Equation \cref{eq:P0_solution} corresponds to the exact solution in the case that $\partial \Omega_R$ is independent of $r_2$, that is, if $f(r_2)$ is a constant. This simplifies the problem significantly since the probability of finding the state associated with the closest molecule of $C$ with a given $r_1$ is independent of the probability of finding the same state with a particular $r_2$. In other words, in this case, Equation \cref{eq:governing_min_eta1_PDE} is separable. However, this simple boundary is equivalent to imposing Smoluchowski's reaction condition on Reaction \cref{eq:example_trimolecular_reaction} and therefore only yields bimolecular mass-action kinetics. To obtain more complex kinetics, in general $\partial \Omega_R$ will need to be spatially dependent, and a separable solution cannot be found. The degree to which the exact density can be approximated by $P_0$ depends on how well $\partial \Omega$ is approximated by the flat boundary assumed when deriving $P_0$. This dependence can be inferred from \cref{fig:flux_gamma_test} which shows that the steady-state reaction rate $K_{\text{FE}}(c)$ converges to the approximate rate $K_1(c)$, which is constructed from $P_0$, as $\Gamma$ reduces since in the limit that $\Gamma \rightarrow 0$ the boundary becomes constant.\par

The approximate steady-state reaction rate $K_1(c)$ neglects $O(\sigmamax^2)$ contributions to the total flux over $\partial \Omega_R$, which amounts to ignoring the $O(\sigmamax^2)$ corrections to $P_0$ and the flux of the probability density over $\partial \Omega_R$ in the $\pos{\hat{r}}{2}$ direction. In essence, we ignore the diffusive and advective flux of the density in the $\pos{\hat{r}}{2}$ direction and we also neglect any influence this diffusive and advective motion has on the probability density itself. As a consequence, compared to $P_h$, $P_0$ underestimates the probability that the closest molecule is close to the origin, see \cref{fig:dif_sigma0.1_gamma1_c10}, while $K_1(c)$ underestimates the contributions of the flux over $\partial \Omega_R$ in both the $\pos{\hat{r}}{1}$ and $\pos{\hat{r}}{2}$ directions, see \cref{fig:flux_sigma_test} and \cref{fig:flux_gamma_test}. Both corrections to $K_1(c)$ are $O(\sigmamax^2)$ but we find that the correction to the flux over $\partial \Omega_R$ in the $\pos{\hat{r}}{1}$ direction is typically larger in magnitude. Although, as shown in \cref{fig:flux_sigma_test}, this is parameter dependent, and it is possible for the correction in the $\pos{\hat{r}}{2}$ direction to be larger.\par

\textcite{our_first_paper} also examined the relative error between $K(c)$ and $K_1(c)$ for the particular choice of $\partial \Omega_R$ given in Equation \cref{eq:exp_boundary}. However, they approximated $K(c)$ using the reaction rate obtained from a particle-based simulation of Reaction \cref{eq:example_trimolecular_reaction}, rather than $K_{\text{FE}}(c)$. As a result, for comparable parameter choices, our values for $\Delta K(c)$ do not agree with those reported in \cite{our_first_paper}. In particular, for $c = 10$ $\sigmamax = 0.1$, $\Gamma = 1$, $\hat{D}_1 = 2$ and $\hat{D}_2 = 1.5$, \cite{our_first_paper} reports a (signed) relative error of less than $-2.5 \times 10^{-2}$ indicating that the simulated rate is less than $K_1(c)$. In comparison, for this choice of parameters, we find that $K_{\text{FE}}(c)$ exceeds $K_1(c)$ and that $\Delta K(c) = 7.6 \times 10^{-3}$. In fact, for all parameter combinations considered, we find that the relative error is strictly positive, suggesting that $K_1(c)$ underestimates the true reaction rate. This disparity arises because the particle-based simulation developed in \cite{our_first_paper} underestimates the actual reaction rate, since it does not exactly replicate the particle dynamics, and so it is possible for reactions to go undetected. Fortunately, such issues can be addressed by using more accurate simulation techniques \cite{egfrd1,egfrd2}, which leaves only the error that arises from using $K_1(c)$ to approximate the reaction rate.\par

In \cref{sec:reactive_boundary_corrections} we demonstrate how to minimise this source of error by altering the reaction boundary $\partial \Omega_R$. To obtain an analytic expression for $\partial \Omega_R$, we use Equation \cref{eq:closest_rate_old} and the inverse Laplace transform to derive a boundary that (approximately) matches the desired reaction rate $K_D(c)$. This assumes that $K_1(c)$ approximates the exact reaction rate, for which $K_{\text{FE}}(c)$ is a proxy, well enough that the corresponding boundary is close enough to the true boundary that reproduces $K_D(c)$ to serve as a good starting point for the optimisation. To retain an analytic expression for $\partial \Omega_R$ we choose to only alter the parameters, $\sigmamax$ and $\Gamma$ in our case, that define the boundary during the optimisation rather than the functional form of the boundary itself. In essence, we are assuming that the true boundary that corresponds to $K_D(c)$ has the same functional form as the initial boundary determined using Equation \cref{eq:closest_rate_old}. This is equivalent to assuming that $K_{\text{FE}}(c)$ for a given $\partial \Omega_R$ has the same dependence on the concentration $c$ as the approximation $K_1(c)$, and only the dependence on the boundary parameters differs between the two rates.\par 

Our approach is simple to implement because it relies on the well-known and widely available Levenberg-Marquardt algorithm, and the optimised boundary is no more difficult to evaluate than the initial boundary. This is particularly important for the efficiency of particle-based simulations, since the value of reaction boundary must be evaluated for each potential reaction event. Moreover, the corrections are independent of the concentration $c$ and so can be precomputed for any given chemical system. However, notice that if the concentration is large, then altering $\partial \Omega_R$ can increase rather than decrease the relative error between $K_{\text{FE}}(c)$ and $K_D(c)$. As shown in \cref{fig:uncorrected_concentration_test}, for the initial boundary, determined by Equation \cref{eq:closest_rate_old}, $K_{\text{FE}}(c)$ will converge to $K_D(c)$ as $c$ increases. This is because the probability density becomes more localised in $r_2$ as c is increased, and hence the portion of $\partial \Omega_R$ over which the flux of the probability density is non-zero continually decreases. Eventually, the relevant segment of $\partial \Omega_R$ becomes small enough that it can be considered essentially constant and, under such conditions, $K_{\text{FE}}(c)$ is equal to $K_1(c)$. That is, for sufficiently large concentrations, Equation \cref{eq:closest_rate_old} yields the correct reaction boundary. Since our approach seeks to minimise the sum of the squares of the residuals across a range of concentrations, the optimisation that will typically improve the error for small concentrations, where the error is initially large, and will potentially make the error worse for large concentrations, where the error is initially small.\par 

It may be possible to achieve better convergence to the desired reaction rate by allowing the functional form of $\partial \Omega_R$ to be varied during the optimisation. An obvious starting point in this direction would be to model the concentration dependence of the error between $K_{\text{FE}}(c)$ and $K_D(c)$ and then attempt to introduce additional terms to the reaction boundary that account for the differing concentration dependence of $K_{\text{FE}}(c)$ and $K_1(c)$. Furthermore, if we abandon our attachment to a simple analytic expression for $\partial \Omega_R$ there is potential to selectively alter portions of the boundary to improve the error for select values of $c$ while minimising the adverse effects for concentrations outside this range. Finally, if the range of $c$ can be estimated in advance, then better results can be obtained by performing the optimisation over this application specific range. However, there exist many biochemical systems in which the concentration of reactants undergo large-scale oscillations \cite{GOLDBETER20132778, SHANKARAN2010650, Kraikivski2021} for which this will not be possible.\par

\section{Conclusion}
We considered the evolution of a system that consists of diffusing molecules from three different chemical species labelled $A$, $B$, and $C$. The molecules react whenever the distance $r_1$ between two molecules of $A$ and $B$ is less than or equal to a function $f(r_2)$ of the distance $r_2$ to the closest molecule of $C$. This proximity-based reaction condition defines a reactive region in the state space of the system on whose surface ($\partial \Omega_R$) triplets of molecules, where each triplet contains one molecule of $A$, $B$ and $C$, are absorbed. We derive a nonlinear partial integro-differential equation, Equation \cref{eq:governing_min_eta1_PDE}, which describes the evolution of the probability density to find the closest molecule, relative to a particular pair of $A$ and $B$ molecules, of $C$ at a given position. The proximity-based reaction condition is imposed as an absorbing boundary condition on the governing equation, see Equation \cref{eq:governing_min_eta1_IB}, and the total flux of the probability density over this boundary corresponds to the reaction rate for the system.\par 

By assuming that $\partial \Omega_R$ extends a small distance $\sigmamax$ in the $\pos{\hat{r}}{1}$ direction, we use singular perturbation theory to derive a leading-order solution ($P_0(r_1,r_2)$) for the steady-state probability density. This solution is then used to construct an approximation, $K_1(c)$ defined in Equation \cref{eq:closest_rate_old}, to the steady-state reaction rate for the system. Both $P_0$ and $K_1(c)$ neglect corrections of $O(\sigmamax^2)$, which are difficult to study analytically, so we construct finite element solutions $P_h$ and $s_h$ for the probability density and the flux of this density respectively, for the specific choice of $\partial \Omega_R$ defined in Equation \cref{eq:exp_boundary}. Using these solutions, we show that $P_0$ underestimates the probability of finding the closest molecule of $C$ near the origin, or at small $r_2$ values as shown in \cref{fig:dif_sigma0.1_gamma1_c10}, while $K_1(c)$ underestimates the reaction rate $K_\text{FE}(c)$ obtained from $s_h$, see \cref{fig:flux_sigma_test}, \cref{fig:flux_gamma_test} and \cref{fig:uncorrected_concentration_test}.\par 

This is the first time such errors have been quantified, and our results elucidate the source of the errors reported in particle-based simulations of an equivalent system \cite{our_first_paper}. In \cref{sec:reactive_boundary_corrections} we demonstrate how to correct for these errors in a manner that can be easily and efficiently incorporated into existing simulations. Moreover, our work provides insight into the dynamics of particles within particle-based simulations and should assist in the continual development of novel proximity-based reaction conditions. Although we have only considered a single functional form for $\partial \Omega_R$, our analysis could be extended to other reactive boundaries. So long as the boundary remains monotonically decreasing in $r_2$ as assumed here, our analysis should extend in a straightforward manner. Violating this assumption complicates matters as the state associated with the closest molecule of $C$ will not necessarily cross the boundary first. If we mandate that reactions only occur when this specific state crosses the boundary, then the dynamics of the system will still be governed by Equation \cref{eq:governing_min_eta1_PDE}, but we would expect the behaviour of the examined errors to change since the normal vector to the boundary will sometimes have a component in the opposite direction of $\pos{\hat{r}}{2}$. In these regions, the advective motion described by the final term in Equation \cref{eq:governing_min_eta1_PDE} will push the relevant state away from the reaction boundary, thus decreasing the reaction rate as opposed to always increasing the rate as is the case for a monotonically decreasing boundary.

\printbibliography

\appendix

\section{The evolution of \texorpdfstring{$\Phi(\boldsymbol{\eta_2},t)$}{Lg}}
\label{sec:AppendixA}
To derive the governing equation for $\PhiFunc$, we consider $\mathcal{L}_2\PhiFunc$ where
\begin{equation}
    \mathcal{L}_2 \equiv \frac{\partial}{\partial t} - \hat{D}_2\hat{\nabla}^2_2,
\end{equation}
is the diffusion operator on the $3$-dimensional space spanned by $\pos{\eta}{2}$, originally defined in Equation \cref{eq:diffusion_eta1}. The time derivative of $\PhiFunc$ is given by
\begin{equation}\label{eq:time_derivative}
    \Phi_t = g\phi_t - c\phi g\int_{V_2}\phi_t\left(\pos{\eta'}{3},t\right) dV_2',
\end{equation}
where we have used the $t$ subscript to denote differentiation with respect to time and defined 
\begin{equation}
    g \equiv \gFunc = \text{exp}\left(-c\int_{V_2}\phiFuncPrime dV_2'\right),
\end{equation}
for notational convenience. We can also take the Laplacian of $\Phi$ with respect to the coordinates of $\pos{\eta}{2}$
\begin{equation}
    \hat{\nabla}^2_2 \Phi = g\hat{\nabla}^2_2\phi + 2\left(\hat{\nabla}_2\phi\right)\cdot\left(\hat{\nabla}_2g\right) + \phi \hat{\nabla}^2_2g,
\end{equation}
which when combined with Equation \cref{eq:time_derivative} yields
\begin{equation}\label{eq:L3_on_Phi}
    \mathcal{L}_2 \Phi = g\phi_t - c\phi g\int_{V_2}\phi_t\left(\pos{\eta'}{3},t\right) dV_2'- \hat{D}_2\left[g\hat{\nabla}^2_2\phi - 2\left(\hat{\nabla}_2\phi\right)\cdot\left(\hat{\nabla}_2g\right) - \phi \hat{\nabla}^2_2g\right].
\end{equation}
\par
Recalling that $\mathcal{L}_2\phi = 0$ from Equation \cref{eq:diffusion_eta1} and applying the Divergence theorem we find
\begin{equation}\label{eq:L3_on_Phi_simplified1}
    \begin{split}
    \mathcal{L}_2 \Phi &= \hat{D}_2\left[-c\phi g\int_{V_2}\hat{\nabla}^2_2\phi\left(\pos{\eta'}{3},t\right) dV_2' - 2\left(\hat{\nabla}_2\phi\right)\cdot\left(\hat{\nabla}_2g\right) - \phi \hat{\nabla}^2_2g\right]\\
    &= - c\hat{D}_2\phi g\oint_{S_2}\left(\hat{\nabla}_2\phi\left(\pos{\eta'}{3},t\right)\right)\cdot \pos{\hat{r}}{2} dA_2' - 2\hat{D}_2\left(\hat{\nabla}_2\phi\right)\cdot\left(\hat{\nabla}_2g\right) - \hat{D}_2\phi \hat{\nabla}^2_2g,
\end{split}
\end{equation}
where $S_2$ is the surface of a sphere of radius $r_2 = ||\pos{\eta}{2}||$, $dA_2'$ is an elemental area on that surface and $\pos{\hat{r}}{2}$ is the unit outward facing normal vector. The diffusion in the $\pos{\eta}{2}$ coordinate is isotropic, and so $\phi$ is independent of orientation. That is, $\phi$ only has radial dependence, so we have
\begin{equation}\label{eq:nabla_g}
\begin{split}
    \hat{\nabla}_2 g &= \hat{\nabla}_2 \text{exp}\left(-c\int_{V_2}\phiFuncPrime dV_2'\right) \\&= -cg\phi\left(\oint_{S_2} dA_2'\right)\pos{\hat{r}}{2}.
\end{split}
\end{equation}
By the same reasoning we are able to take the integrand outside the integral in Equation \cref{eq:L3_on_Phi_simplified1} and by substituting in Equation \cref{eq:nabla_g} we obtain
\begin{equation}
    \begin{split}
    \mathcal{L}_2 \Phi &= -\hat{D}_2 \parentheses{[}{\parentheses{(}{\hat{\nabla}_2 \phi}{)}\cdot \parentheses{(}{\hat{\nabla}_2g}{)} + \phi \hat{\nabla}^2_2g}{]}\\
    &= -\hat{D}_2\hat{\nabla}_2\cdot \parentheses{(}{\phi\hat{\nabla}_2g}{)}\\
    &= \hat{D}_2\hat{\nabla}_2 \cdot \parentheses{(}{cg\phi^2\left(\oint_{S_2} dA_2'\right)\pos{\hat{r}}{2}}{)}\\
    &= \hat{D}_2\hat{\nabla}_2 \cdot \parentheses{(}{4\pi r_2^2 c\phi\Phi\pos{\hat{r}}{2}}{)}.
\end{split}
\end{equation}
That is,
\begin{equation}
    \frac{\partial \PhiFunc}{\partial t} = \hat{D}_2\hat{\nabla}^2_2\PhiFunc + \hat{D}_2\hat{\nabla}_2 \cdot \parentheses{(}{4\pi r_2^2 c\phi\PhiFunc\pos{\hat{r}}{2}}{)},
\end{equation}
as stated in Equation \cref{eq:Phi_evolution}.

\end{document}